\documentclass[a4paper,12pt]{amsart}
\usepackage{amsfonts}
\usepackage{amssymb}
\usepackage{ifthen}
\usepackage{graphicx}
\nonstopmode \numberwithin{equation}{section}
\usepackage[usenames]{color}
\usepackage{enumerate}
\usepackage{graphicx}
\usepackage{epstopdf}
\numberwithin{equation}{section}

\theoremstyle{definition}

\newtheorem{examples}[equation]{Examples}

\newtheorem{thm}{Theorem}[section]
\newtheorem{lem}[thm]{Lemma}
\newtheorem{cor}[thm]{Corollary}

\newtheorem{rem}[thm]{Remark}

\newtheorem{cl}{Claim}[section]
\newtheorem{ca}{Case}[section]
\newtheorem{sca}{Subcase}[section]
\newtheorem{scl}[section]{Subclaim}
\newtheorem{conj}[equation]{Conjecture}

\theoremstyle{definition}
\newtheorem{defn}[thm]{Definition}
\newtheorem{op}[equation]{Open Problem}
\newtheorem{ques}[thm]{Question}
\newtheorem{exam}[equation]{Example}

\newcounter {own}
\def\theown {\thesection       .\arabic{own}}

\newenvironment{pf}[1][]{%
 \vskip 3mm
 \noindent
 \ifthenelse{\equal{#1}{}}%
  {{\slshape Proof. }}%
  {{\slshape #1.} }%
 }%
{\qed\bigskip}

\newcounter{alphabet}

\newcommand{\IR}{{\mathbb R}}

\newcommand{\diam}{{\operatorname{diam}}}

\newcommand{\dist}{{\operatorname{dist}}}

\def\be{\begin{equation}}
\def\ee{\end{equation}}

\newcommand{\bee}{\begin{enumerate}}
\newcommand{\eee}{\end{enumerate}}

\newcommand{\blem}{\begin{lem}}
\newcommand{\elem}{\end{lem}}
\newcommand{\bthm}{\begin{thm}}
\newcommand{\ethm}{\end{thm}}
\newcommand{\bcor}{\begin{cor}}
\newcommand{\ecor}{\end{cor}}
\newcommand{\beg}{\begin{exam}}
\newcommand{\eeg}{\end{exam}}
\newcommand{\begs}{\begin{examples}}
\newcommand{\eegs}{\end{examples}}
\newcommand{\bdefe}{\begin{defn}}
\newcommand{\edefe}{\end{defn}}
\newcommand{\bprob}{\begin{prob}}
\newcommand{\eprob}{\end{prob}}
\newcommand{\bques}{\begin{ques}}
\newcommand{\eques}{\end{ques}}
\newcommand{\bei}{\begin{itemize}}
\newcommand{\eei}{\end{itemize}}
\newcommand{\bcon}{\begin{conj}}
\newcommand{\econ}{\end{conj}}
\newcommand{\bop}{\begin{op}}
\newcommand{\eop}{\end{op}}

\newcommand{\bca}{\begin{ca}}
\newcommand{\eca}{\end{ca}}
\newcommand{\bsca}{\begin{sca}}
\newcommand{\esca}{\end{sca}}

\newcommand{\bcl}{\begin{cl}}
\newcommand{\ecl}{\end{cl}}

\newcommand{\bscl}{\begin{scl}}
\newcommand{\escl}{\end{scl}}

\newcommand{\bcons}{\begin{conjs}}
\newcommand{\econs}{\end{conjs}}
\newcommand{\bprop}{\begin{propo}}
\newcommand{\eprop}{\end{propo}}
\newcommand{\br}{\begin{rem}}
\newcommand{\er}{\end{rem}}
\newcommand{\brs}{\begin{rems}}
\newcommand{\ers}{\end{rems}}
\newcommand{\bo}{\begin{obser}}
\newcommand{\eo}{\end{obser}}
\newcommand{\bos}{\begin{obsers}}
\newcommand{\eos}{\end{obsers}}
\newcommand{\bpf}{\begin{pf}}
\newcommand{\epf}{\end{pf}}
\newcommand{\ba}{\begin{array}}
\newcommand{\ea}{\end{array}}
\newcommand{\beq}{\begin{eqnarray}}
\newcommand{\beqq}{\begin{eqnarray*}}
\newcommand{\eeq}{\end{eqnarray}}
\newcommand{\eeqq}{\end{eqnarray*}}

\newcounter{minutes}
\divide\time by 60
\newcounter{hours}
\multiply\time by 60 \addtocounter{minutes}{-\time}
\begin{document}

\bibliographystyle{amsplain}
\title{Pommerenke's theorem on Gromov hyperbolic domains}

\author{Qingshan Zhou$^{~\mathbf{*}}$}
\address{Qingshan Zhou, School of Mathematics and Big Data, Foshan university,  Foshan, Guangdong 528000, People's Republic
of China} \email{qszhou1989@163.com; q476308142@qq.com}

\author{Antti Rasila}
\address{Antti Rasila, Technion -- Israel Institute of Technology, Guangdong Technion, Shantou, Guangdong 515063, People's Republic
of China} \email{antti.rasila@gtiit.edu.cn; antti.rasila@iki.fi}

\author{Tiantian Guan}
\address{Tiantian Guan, Department of Mathematics,
Shantou University, Shantou, Guangdong 515063, People's Republic of
China} \email{ttguan93@163.com}

%%%%%%%% BEGIN TIMESTAMP
\def\thefootnote{}
\footnotetext{ \texttt{\tiny File:~\jobname .tex,
          printed: \number\year-\number\month-\number\day,
          \thehours.\ifnum\theminutes<10{0}\fi\theminutes}
} \makeatletter\def\thefootnote{\@arabic\c@footnote}\makeatother
%%%%%%%% END TIMESTAMP

\date{}
 \keywords{
Gromov hyperbolic domain, quasihyperbolic geodesic, Gehring-Hayman property, uniform domain, quasisymmetric, Busemann function. \\ ${}^{\mathbf{*}}$ Corresponding author
}

\begin{abstract}
We establish a version of a classical theorem of Pommerenke, which is a diameter version of the Gehring-Hayman inequality on Gromov hyperbolic domains of $\mathbb{R}^n$. Two applications are given. Firstly, we generalize Ostrowski's Faltensatz to quasihyperbolic geodesics of Gromov hyperbolic domains. Secondly, we prove that unbounded uniform domains can be characterized in the terms of Gromov hyperbolicity and a naturally quasisymmetric correspondence on the boundary, where the Gromov boundary is equipped with a Hamenst\"adt metric (defined by using a Busemann function).
\end{abstract}

\thanks{Qingshan Zhou was supported by NSF of China (No. 11901090). Antti Rasila was supported by NSF of China (No. 11971124) and NSF of the Guangdong Province (No.  2021A1515010326).}

\maketitle{} \pagestyle{myheadings} \markboth{Qingshan Zhou, Antti Rasila, and Tiantian Guan}{Pommerenke's theorem on Gromov hyperbolic domains}

\section{Introduction and main results}\label{sec-1}
This paper deals with geometric properties of quasihyperbolic geodesics on Gromov hyperbolic domains of Euclidean spaces $\mathbb{R}^n$ ($n\geq 2$). Let $G\subsetneq \mathbb{R}^n$ be a domain (an open and connected set). The quasihyperbolic metric $k=k_G$ of $G$ was introduced by Gehring and Osgood in \cite{GO}. It is well known that $(G, k)$ is complete, proper, and geodesic as a metric space, see \cite[Proposition 2.8]{BHK}. We say that $G$ is a {\it Gromov hyperbolic domain} if $(G, k)$ is a $\delta$-hyperbolic space in the sense of Gromov (see Subsection \ref{sec-g}), for some constant $\delta\geq 0$. It is known that this class of domains include:
\begin{enumerate}
\item simply connected planar domains with nonempty boundaries;
\item bounded convex domains of $\mathbb{R}^n$, for example, the unit balls $\mathbb{B}^n$;
\item uniform domains which were introduced by Martio and Sarvas \cite{MS78} and Jones \cite{Jo80};
\item the image domains of uniform domains and $\mathbb{B}^n$ under a quasiconformal mapping. Note that Gromov hyperbolicity of domains in $\mathbb{R}^n$ is preserved under quasiconformal mappings, see \cite[Page 5]{BHK}.
\end{enumerate}

In \cite[Section 7]{BHK}, Bonk, Heinonen, and Koskela conjectured that there is a characterization of Gromov hyperbolic domains in terms of the Euclidean geometry. This conjecture has been proved true by Balogh and Buckley \cite{BB03} in the setting of metric measure spaces with a bounded geometry.
\begin{thm}$($\cite[Theorem 0.1]{BB03}$)$\label{BB03-thm0.1}
Let $G\subsetneq \mathbb{R}^n$ be a domain. Then $G$ is Gromov $\delta$-hyperbolic if and only if the following statements hold:
\begin{enumerate}
\item Gehring-Hayman property: there is a constant $C_{gh}\geq 1$ such that for any quasihyperbolic geodesic $[x, y]_k$ and for every curve $\gamma$ in $G$ with end points $x$ and $y$, we have that
    $$\ell([x, y]_k)\leq C_{gh}\ell(\gamma),$$
    where $\ell(\gamma)$ denotes the Euclidean length of $\gamma$;
\item Separation property: there is a constant $C_{sp}\geq 1$ such that for any quasihyperbolic geodesic $[x, y]_k$, for any point $z\in [x, y]_k$,  and for every curve $\gamma$ in $G$ with end points $x$ and $y$, then
    $$B_{\ell_{G}}(z, C_{sp}d(z))\cap \gamma \neq \emptyset,$$
    where $d(z)=\dist(z, \partial G)$, and $\ell_G$ is the length metric of $G$ associated with the $n$-Euclidean metric $|\cdot|$, and $B_{\ell_G}\big(z, C_{sp}d(z)\big)=\{x\in G\,|\,\ell_G(x, z)<C_{sp}d(z)\}$,
\end{enumerate}
where $\delta$, $C_{gh}$ and $C_{sp}$ depend only on each other and $n$.
\end{thm}

The first property was observed by Gehring and Hayman \cite{GH}. It asserts  that the hyperbolic geodesic as a near-minimal length among all curves connecting its endpoints
in a simply connected planar domain. There are several generalizations of this result, see \cite{BHK, KL} for more information.

 In \cite{P64}, Pommerenke proved a diameter version of the Gehring-Hayman theorem for hyperbolic geodesics. By studying distortion properties of arcs, Heinonen and N\"akki \cite{HN} found two extensions of Pommerenke's theorem. In the dimension $n=2$, the quasihyperbolic geodesic essentially has the smallest diameter among all curves with the same end points in a simply connected domain, and in higher dimensions $n\geq 3$, they extended the result to a domain which is quasiconformally equivalent to a uniform domain, see \cite[Theorems 6.2 and 6.8]{HN}.

It is reasonable to ask whether the theorem of Pommerenke holds for quasihyperbolic geodesics of Gromov hyperbolic domains in $\mathbb{R}^n$. We show the following as our main result:
\begin{thm}\label{main thm-1}
Let $G\subsetneq \mathbb{R}^n$ be a $\delta$-hyperbolic domain. Then for any quasihyperbolic geodesic $[x, y]_k$ and for each curve $\gamma$  in $G$ with end points $x$ and $y$, there is a constant $C_{po}=C_{po}(n,\delta)$ such that
$$\diam([x, y]_k)\leq C_{po}\diam(\gamma).$$
If in addition, $G$ is a domain of $C_{bt}$-bounded turning, then
$$\diam([x, y]_k)\leq C_{bt}C_{po}|x-y|.$$
\end{thm}

As an application of Theorem \ref{main thm-1}, we show Ostrowski's Faltensatz for quasihyperbolic geodesics of Gromov hyperbolic domains of higher dimensional Euclidean spaces. The classical Ostrowski's Faltensatz  \cite{O} states that:
\begin{thm}$($\cite[Page 318]{P75}$)$\label{O}
Let $G\subsetneq \mathbb{R}^2$ be a simply connected domain, let $\Sigma$
be a cross-cut of $G$ dividing it into two subdomains $G_1$ and $G_2$, and let $L$ be
a hyperbolic geodesic in $G$ whose end points lie in $G_1$. Then every point
$x\in L\cap G_2$ satisfies
$$\dist(x, \Sigma) \leq A \min\{d(x), \diam(\Sigma)\},$$
where $A>0$ is an absolute constant.
\end{thm}
It was showed by Heinonen and N\"akki that the version of Theorem \ref{O} is valid for quasihyperbolic geodesics both in simply connected planar domains and domains in $\mathbb{R}^n$ which are quasiconformally equivalent to uniform domains, see \cite[Theorem 7.3 and Corollary 7.4]{HN}.

Using Theorem \ref{main thm-1} we prove the following  analogy of Ostrowski's Faltensatz on Gromov hyperbolic domains:
\begin{thm}\label{thm-4}
Let $G\subsetneq \mathbb{R}^n$ be a $\delta$-hyperbolic domain, let $\Sigma$ be a cross-section of $G$ so that $G_1$ and $G_2$ are the components of $G\setminus \Sigma$, and let $L$ be a quasihyperbolic geodesic in $G$ whose end points lie in $G_1$. Then there is a constant $A=A(n, \delta)$ such that every point $x\in L\cap G_2$ can be joined to $\Sigma$ by an arc $\alpha$ in $G$ satisfying
$$\ell(\alpha)\leq A \min\{d(x), \diam(\Sigma)\}.$$
In particular,
$$\dist(x, \Sigma)\leq A\min\{d(x), \diam(\Sigma)\}.$$
\end{thm}
\br
Because simply connected planar domains and domains of $\mathbb{R}^n$, which are quasiconformally equivalent to uniform domains, are Gromov hyperbolic, Theorem \ref{thm-4} is an improvement of \cite[Theorem 7.3 and Corollary 7.4]{HN}. Furthermore, we see from the proof that the curve $\alpha$ in the theorem can be chosen by a quasihyperbolic geodesic of $G$ between $x$ and the cross-section $\Sigma$.
\er

\subsection{Gromov hyperbolicity and unbounded uniform domains}
Next we investigate the relationship between Gromov hyperbolic domains and unbounded uniform domains in $\mathbb{R}^n$. It was shown in \cite[Theorem 1.11]{BHK} that a bounded domain in $\mathbb{R}^n$ is uniform if and only if it is both Gromov hyperbolic and its Euclidean boundary is naturally quasisymmetrically equivalent to the Gromov boundary. It is natural to seek an analogue of this result for unbounded domains.

The challenge is that the Gromov boundary endowed with a visual metric (based at a point in the interior) is a bounded metric space. However, quasisymmetric mappings map unbounded sets onto unbounded sets. To overcome this obstacle, we introduce another class of Hamenst\"adt metrics on punctured Gromov boundary via Busemann functions, see Section \ref{sec-3} for more information.

By Theorem \ref{main thm-1}, we obtain the following characterization of unbounded uniform domains in terms of Gromov hyperbolicity:

\begin{thm}\label{main thm-3}
Let $G\subsetneq \mathbb{R}^n$ be an unbounded domain. Then $G$ is $A$-uniform if and only if $G$ is $\delta$-hyperbolic and there is a natural $\eta$-quasisymmetric identification
$$\varphi:\, (\partial G, |\cdot|)\to (\partial_\infty G\setminus\{\xi_0\}, d_{b, \varepsilon_0}),$$
where $\partial_\infty G$ is the Gromov boundary of $(G, k)$, $\varphi(\infty)=\xi_0$, $b\in \mathcal{B}(\xi_0)$ is a Busemann function, and $d_{b, \varepsilon_0}$ is a Hamenst\"adt metric based at $b$ with parameter $\varepsilon_0=\varepsilon_0(\delta) >0$. The constants $A$ and $\delta$, $\eta$ depend on each other and $n$.
\end{thm}
\br
Our strategy in proving Theorem \ref{main thm-3} is the following:
\begin{enumerate}
\item For the necessity, it was observed in \cite[Proposition 3.12]{BHK} that there is a natural identification $\varphi:\, \partial_\infty G\to \partial G\cup \{\infty\}$ between Gromov boundary and the extended Euclidean boundary of $G$ with $\varphi(\xi_0)=\infty$ for some $\xi_0\in \partial_\infty G$. We endow the punctured Gromov boundary $\partial_\infty G\setminus\{\xi_0\}$ with a Hamenst\"adt metric based at a Busemann function $b\in \mathcal{B}(\xi_0)$. Then the assertion follows from a careful computation.
    \item We show the sufficiency part in three steps. First, we consider the conformal deformation of $G$ which was recently introduced in \cite{Z20} and defined by using a density $\rho_\varepsilon(x)=e^{-\varepsilon b(x)}$ for some $\varepsilon>0$ and a Busemann function $b\in \mathcal{B}(\xi_0)$, where $x\in G$. The resulting length space $G_\varepsilon=(G, d_\varepsilon)$ is an unbounded uniform space such that there is a natural quasisymmetric correspondence between $(\partial_\infty G\setminus \{\xi_0\}, d_{b, \varepsilon_0})$ and the metric boundary $\partial G_\varepsilon$. Second, we prove that the domain $G$ is linearly locally connected (LLC)  via the preservation of  modulus of curve families under the conformal deformation.  Finally, it is not difficult to show from
        Theorem \ref{main thm-1} that each quasihyperbolic geodesic in a LLC Gromov hyperbolic domain is a uniform arc. Note that this statement extends \cite[Proposition 7.12]{BHK} to the unbounded case.
\end{enumerate}
\er

The rest of this paper is organized as follows. In Section \ref{sec-2}, we recall some definitions and preliminary results. Section \ref{sec-3} concerns the geometric properties of Gromov hyperbolic domains. The proofs of Theorem \ref{main thm-1} and Theorem \ref{thm-4} are given in Section \ref{sec-4}. Finally, Section 5 is devoted to the proof of Theorem \ref{main thm-3}.

\section{Preliminaries}\label{sec-2}

%\subsection{Notation} %Following \cite{BuSc}, for $a,b,c\in \mathbb{R}$ with $c\geq 0$, it is convenient to write
%$$a{\color{red}\doteq} b\;\;\mbox{up to an error $\leq c$ or}\;\; a{\color{red}\doteq}_c b\;\mbox{instead %of}\;\;|a-b|\leq c.$$
%Also, we extend our argument about ${\color{red}\doteq}$ as follows. For any sequences $\{a_i\},\{b_i\}\subset \mathbb{R}$, we write
%$$\{a_i\}_i {\color{red}\doteq} \{b_i\}_i\;\;\mbox{up to an error $\leq c$ or}\;\; \{a_i\}_i {\color{red}\doteq}_c \{b_i\}_i$$
%if $\limsup\limits_{i\to \infty} |a_i-b_i|\leq c$. Similarly, we denote
%$$\{a_i\}_i\leq \{b_i\}_i+c,\;\;\mbox{if}\;\limsup_{i\to \infty} a_i\leq \liminf_{i\to \infty} b_i+c$$
%and
%$$\{a_i\}_i\geq \{b_i\}_i+c,\;\;\mbox{if}\;\liminf_{i\to \infty} a_i\geq \limsup_{i\to \infty} b_i+c.$$

%Moreover, let $a,b>0$ and $c\geq 1$. We use the notation $a{\color{red}\asymp} b$ up to multiplicative error $\leq c$ or $a{\color{red}\asymp}_c b$ instead of $1/c\leq a/b\leq c$.

%For real numbers $s,t$, we set
%$$s{\color{red} \wedge} t=\min\{s,t\}\;\;\;\;\mbox{and}\;\;\;\; s{\color{red}\vee} t=\max\{s,t\}.$$

\subsection{Metric geometry}
Let $(X,d)$ be a metric space. It is {\it non-complete} if its boundary $\partial X= \overline{X}\setminus X\not=\emptyset$, where $\overline{X}$ denotes the metric completion of $X$. Also, $B(x,r)$  and $\overline{B}(x,r)$ are the open ball and closed ball (of radius $r$ centered at the point $x$) in $(X,d)$, respectively. The metric sphere $S(x, r)=\{y\in X\,|\,d(y, x)=r\}$. $X$ is called {\it proper} if its closed balls are compact. The diameter of a set $E\subset X$ denotes by $\diam(E)$. For all sets $E_1$ and $E_2$ in $(X, d)$, $\dist_d(E_1, E_2)$ means the distance between $E_1$ and $E_2$.

A {\it curve} in $X$ means a continuous map $\gamma:\; I\to X$ from an interval $I\subset \IR$ to $X$. If $\gamma$ is an embedding of $I$, it is called an {\it arc}. We also denote the image set $\gamma(I)$ of $\gamma$ by $\gamma$. The {\it length} $\ell(\gamma)$ of $\gamma$ with respect to the metric $d$ is defined in an obvious way. Here the parameter interval $I$ is allowed to be open or half-open.
We also denote the subarc of $\gamma$ by $\gamma[x,y]$ with end points $x$ and $y$ in $\gamma$ and $\gamma(x, y)=\gamma[x, y]\setminus\{x, y\}$.  Also, $X$ is called {\it rectifiably connected} if every pair of points in $X$ can be joined with a curve $\gamma$ in $X$ with $\ell(\gamma)<\infty$.

A {\it geodesic} $\gamma$ joining $x$ to $y$ in $X$ is a map $\gamma:I=[0,l]\to X$ from an interval $I$ to $X$ such that $\gamma(0)=x$, $\gamma(l)=y$ and $d(\gamma(t),\gamma(t'))=|t-t'|$ for all $t,t'\in I$. A metric space $X$ is said to be {\it geodesic} if every pair of points can be joined by a geodesic arc.

\subsection{Quasihyperbolic metric and uniform spaces}
Let $(X, d)$ be a non-complete, locally compact, and rectifiably connected metric space. The {\it quasihyperbolic metric} $k$ in $X$ is defined by
$$k(x, y)=\inf \Big\{\int_{\gamma} \frac{1}{d(z)}\;ds\Big\},$$
where the infimum is taken over all rectifiable curves $\gamma$ in $X$ with the end points $x$ and $y$, $d(z):=\dist(z,\partial X)$ and $ds$ denotes the arc-length element with respect to the metric $d$. We need the following fact concerning the quasihyperbolic metric:
\begin{lem}$($\cite[Theorem 3.9]{Vai-4}$)$\label{Vai-4 thm3.9}
Let $G\subsetneq \mathbb{R}^n$ be a domain, and let $x,$ $y\in G$ such that $|x-y|\leq 1/2\,d(x)$ or $k(x, y)\leq 1$. Then,
$$\frac{1}{2}\frac{|x-y|}{d(x)}\leq k(x, y)\leq 2\frac{|x-y|}{d(x)}.$$
\end{lem}

\bdefe
A non-complete, locally compact, and rectifiably connected metric space $(X, d)$ is called {\it $A$-uniform} if each pair of points $x, y\,\in X$ can be joined by an $A$-uniform curve $\gamma$. A curve $\gamma$ with end points $x$ and $y$ is said to be an {\it $A$-uniform curve} if:
\begin{enumerate}
\item (Quasiconvex condition) $\ell(\gamma)\leq Ad(x, y)$;
\item (Double cone condition) $\min\{\ell(\alpha[x,z]), \ell(\alpha[z,y])\}\leq A d(z)$ for all $z\in \gamma$.
\end{enumerate}
If the lengths in the definition are replaced by diameters, then it is called an {\it $A$-diameter uniform curve}. Moreover, $X$ is called {\it diameter uniform} if there is a constant $A>0$ such that each pair of points can be joined by an $A$-diameter uniform curve.
\edefe

%\begin{thm}\label{Thm-3}
%Let $X$ be a locally compact, rectifiably connected and noncomplete metric spaces. \begin{enumerate}
%\item\cite[Proposition 2.8]{BHK} If the identification $(X, d)\to (X, l)$ is  a homeomorphism. Then $(X, k_X)$ is proper and geodesic;
%\item\cite[Theorem 2.10]{BHK}\label{e3-2} A quasihyperbolic geodesic in an $A$-uniform space is a $A_1$-uniform curve with $A_1$ depending only on $A$.
%\end{enumerate}
%\end{thm}

\begin{thm}$($\cite[(2.4) and (2.16)]{BHK}$)$\label{BHK 2.4}
Let $(X, d)$ be an $A$-uniform space. Then
 $$\log\bigg(1+\frac{d(x, y)}{\min\{d(x), d(y)\}}\bigg)\leq k(x, y)\leq 4A^2\log\bigg(1+\frac{d(x, y)}{\min\{d(x), d(y)\}}\bigg)$$
 for each pair of points $x, y\in X$.
\end{thm}

\subsection{Modulus and Loewner spaces}
In this part, we assume that $(X, d)$ is a rectifiably connected metric space, and $\mu$ is a Borel measure on $X$.
\bdefe
Given $Q> 1$, we say that
$X$ is {\it Ahlfors $Q$-regular} if there exists a constant $C>0$ such that
for each $x\in X$ and $0<r\leq \diam  (X)$,
$$C^{-1}r^Q\leq \mu(B(x,r))\leq Cr^Q.$$
%It is well known that $\mathbb{R}^n$ is $n$-Ahlfors regular with the $n$-dimensional Lebesgue measure.

If $X$ is non-complete, then $X$ is said to be {\it locally Ahlfors $Q$-regular} provided there is a number $\lambda\in (0,1/2)$ such that the Ahlfors $Q$-regularity holds for all  $x\in X$ and all $r\in (0,\lambda d(x))$.
\edefe

\bdefe\label{de-1}
Let $Q>1$. We define the {\it $Q$-modulus} of a family $\Gamma$ of curves in a metric measure space $(X, d,\mu)$ by
\be\label{eq2-0}\mbox{mod}_Q\;\Gamma=\inf \int_X \rho^Q \; d\mu,\ee
where the infimum is taken over all Borel functions $\rho:X\to [0,\infty]$ satisfying
$$\int_\gamma \rho \; ds\geq 1,$$
 for all locally rectifiable curves $\gamma\in \Gamma$. Then the $Q$-modulus  of a pair of disjoint compact sets $E,F\subset X$ is
 $$\mbox{mod}_Q\;(E,F;X)=\mbox{mod}_Q\;\Gamma(E,F;X),$$
 where $\Gamma(E,F;X)$ is the family of all curves in $X$ joining the sets $E$ and $F$.
\edefe

The notion of {\it a Loewner space} was introduced by Heinonen and Koskela \cite{HK} in their study of the equivalence of quasiconformality and quasisymmetry. Note that a  metric measure space $(X, d,\mu)$ is called $Q$-{\it Loewner}, with $Q>1$, provided that the Loewner control function
$$\phi(t)=\inf \Big\{\mbox{mod}_Q\;(E,F;X)\,|\, \Delta(E,F)\leq t\Big\}$$
is strictly positive for all $t>0$; here $E$ and $F$ are non-degenerate disjoint continua (connected and compact sets) in $X$, and
$$\Delta(E, F) = \frac{\dist(E, F)}{\min\{\diam (E), \diam (F)\}}.$$
It is known that $\mathbb{R}^n$ and the unit balls $\mathbb{B}^n$ are $n$-Loewner with the control function depending only on $n$ (see \cite[Section 6]{HK}).

Following [3], let $(X, d, \mu)$ be a locally compact and non-complete metric measure space, and let $Q > 1.$ We say that $X$ is {\it  locally $Q$-Loewner} if there exist $\kappa\geq1$, $\lambda_0\in(0,1]$ and a decreasing function $\psi: (0,\infty) \to (0,\infty)$ such that for all $0<\lambda\leq \lambda_0$, for all $x\in X$, and for all non-degenerate disjoint continua $E$ and $F$ in $B(x, \lambda d(x))$ with $\Delta(E,F)\leq t$,
$$\mbox{mod}_Q\Big(E,F;B\big(x, \kappa\lambda d(x)\big) \Big) \geq \psi(t).$$

%Then we record certain auxiliary results concerning the curve family modulus as follows.

%\begin{thm}\cite[Lemma 3.14]{HK}\label{z5}
% Suppose that a metric measure space $(X, d,\mu)$ is a Loewner space of Hausdorff dimension $Q>1$ satisfying the upper regular. Let $y\in X$ and let $0<2r<R$. Then $$\mbox{mod}_Q\Big(\overline{B}(y,r), X\setminus B(y,R)\Big)\leq C_\mu\Big(\log \frac{R}{r}\Big)^{1-Q}.$$
%\end{thm}

\subsection{Mappings on metric spaces}
%%%%%%%%%%%%%%%%%%%

Following \cite{BHK,HRWZ}, a homeomorphism $f:$ $(X,d)\to (X',d')$ between two metric spaces is said to be {\it $\eta$-quasisymmetric} if there is a homeomorphism $\eta : [0,\infty) \to [0,\infty)$ such that
$$ d(x,a)\leq t d(x,b)\;\; \mbox{implies}\;\;   d'(f(x),f(a)) \leq \eta(t) d'(f(x),f(b))$$
for each $t>0$ and for each triplet $x,$ $a$, $b$ of points in $X$. A homeomorphism $f:$ $(X,d)\to (X',d')$ between two non-complete metric spaces is said to be {\it $q$-locally $\eta$-quasisymmetric} if there is a constant $q\in (0,1)$ and a homeomorphism $\eta : [0,\infty) \to [0,\infty)$ such that $f$ is $\eta$-quasisymmetric in $B(x, qd(x))$ for all $x\in X$.

\subsection{Gromov hyperbolicity}\label{sec-g}
Let $(X,d)$ be a metric space. Fix a base point $o\in X$. For any two points $x,y\in X$,  we define
$$(x|y)_o=\frac{1}{2}\big(d(x,o)+d(y,o)-d(x,y)\big).$$
This number is called the  {\it Gromov product} of $x$ and $y$ with respect to $o$.

We say that $(X, d)$ is a {\it Gromov hyperbolic space}, if there is a constant $\delta\geq 0$ such that
$$(x|y)_o\geq \min\{(x|z)_o ,(z|y)_o\}-\delta,$$
for all $x,y,z,o\in X$.

Suppose that $(X, d)$ is a Gromov $\delta$-hyperbolic metric space for some constant $\delta\geq 0$. A sequence $\{x_i\}$ in $X$ is called a {\it Gromov sequence} if $(x_i|x_j)_o\rightarrow \infty$ as $i,$ $j\rightarrow \infty.$ Two such Gromov sequences $\{x_i\}$ and $\{y_j\}$ are said to be {\it equivalent} if $(x_i|y_i)_o\rightarrow \infty$ as $i\to\infty$. The {\it Gromov boundary} $\partial_\infty X$ of $X$ is defined to be the set of all equivalence classes of Gromov sequences.

 For $x\in X$ and $\xi\in \partial_\infty X$, the Gromov product $(x|\xi)_o$ of $x$ and $\xi$ is defined by
$$(x|\xi)_o= \inf \big\{ \liminf_{i\rightarrow \infty}(x|y_i)_o\;\big |\; \{y_i\}\in \xi\big\}.$$
For $\xi,$ $\zeta\in \partial_\infty X$, the Gromov product $(\xi|\zeta)_o$ of $\xi$ and $\zeta$ is defined by
$$(\xi|\zeta)_o= \inf \big\{ \liminf_{i\rightarrow \infty}(x_i|y_i)_o\;\big |\; \{x_i\}\in \xi\;\;{\rm and}\;\; \{y_i\}\in \zeta\big\}.$$

Also, we need auxiliary results on the Gromov products.

\begin{lem}$($\cite[Standard estimate 2.33]{Vai-0}$)$\label{2.33}
Suppose that $X$ is a geodesic $\delta$-hyperbolic metric space, and $o\in X$. For any geodesic $\gamma$ with end points $x$ and $y$,
$$\dist(o, \gamma)-2\delta \leq (x|y)_o\leq \dist(o, \gamma).$$
\end{lem}

\begin{lem}\label{z0}$($\cite[Lemma $5.11$]{Vai-0}$)$
Let $o,z\in X$, let $X$ be a $\delta$-hyperbolic space, and let $\xi,\xi'\in\partial_\infty X$. Then for any sequences $\{y_i\}\in \xi$, $\{y_i'\}\in \xi'$, we have
\begin{enumerate}
\item  $(z|\xi)_o\leq  \liminf\limits_{i\rightarrow \infty} (z|y_i)_o \leq  \limsup\limits_{i\rightarrow \infty} (z|y_i)_o\leq (z|\xi)_o+\delta;$
\item  $(\xi|\xi')_o\leq  \liminf\limits_{i\rightarrow \infty} (y_i|y_i')_o \leq  \limsup\limits_{i\rightarrow \infty} (y_i|y_i')_o\leq (\xi|\xi')_o+2\delta.$
\end{enumerate}
\end{lem}

%Let $X$ be a $\delta$-hyperbolic space and $o\in X$ be given. For $0<\varepsilon<\min\{1,\frac{1}{5\delta}\}$, we define
%$$\rho_{o,\varepsilon}(\xi,\zeta)=e^{-\varepsilon(\xi|\zeta)_o}$$
%for all $\xi,\zeta$ in the Gromov boundary $\partial_\infty X$ of $X$ with convention $e^{-\infty}=0$.

%We now define
%$$d_{o,\varepsilon} (\xi,\zeta):=\inf\Big\{\sum_{i=1}^{n} \rho_{o,\varepsilon} (\xi_{i-1},\xi_i)\;\big |\;n\geq 1,\xi=\xi_0,\xi_1,\ldots,\xi_n=\zeta\in \partial_\infty X\Big\}.$$
%Then $(\partial_\infty X,d_{o,\varepsilon})$ is a metric space with
%\be\label{b-3} \frac{1}{2}\rho_{o,\varepsilon}\leq d_{o,\varepsilon}\leq \rho_{o,\varepsilon},\ee
%and we call $d_{o,\varepsilon}$ the {\it visual metric} or {\it Bourdon metric} of $\partial_\infty X$ with base point $o\in X$ and the parameter $\varepsilon$.

Finally, we conclude this subsection by giving a definition which was introduced by V\"{a}is\"{a}l\"{a} in \cite{Vai}.
\bdefe
Let $X$ be a proper and geodesic $\delta$-hyperbolic space, and let $K\geq 0$.
Let $\xi\in\partial_\infty X$. We say that $X$ is {\it $K$-roughly starlike} with respect to $\xi$ if for each $x\in X$, there is $\eta\in\partial_\infty X$ and a geodesic $\gamma$ between $\xi$ and $\eta$ such that  $\dist(x,\gamma)\leq K.$
\edefe

\subsection{Busemann functions and the Hamenst\"adt metric}
 Let $(X, d)$ be a $\delta$-hyperbolic space, and let $\xi\in\partial_\infty X$, and fix $o\in X$. We say that $b:X\to \mathbb{R}$ is a {\it Busemann function} based at $\xi$, denoted by $b\in \mathcal{B}(\xi)$, if for all $x\in X$, we have
$$b(x):=b_{\xi,o}(x):=b_\xi(x,o)=(o|\xi)_x-(x|\xi)_o.$$

By \cite[Proposition 3.1.5(1)]{BuSc}, for all $x,y\in X$ we have
\be\label{b-1} |b(x)-b(y)|\leq d(x,y)+10\delta.\ee
It follows from \cite[Lemma $3.1.1$]{BuSc} that
\be\label{zz0}
b(x)-2\delta\leq \limsup\limits_{i\to \infty} \big(d(x,z_i)-d(o,z_i)\big)  \leq b(x)+2\delta\ee
for every Gromov sequence $\{z_i\}\in\xi$.

We next define the Gromov product of $x,y\in X$ based at $b$ by
$$(x|y)_b=\frac{1}{2}(b(x)+b(y)-d(x,y)).$$
Moreover, by \cite[(3.2) and Example 3.2.1]{BuSc}, we see that
\be\label{zz0.1}
(x|y)_b-10\delta\leq (x|y)_o-(x|\xi)_o-(y|\xi)_o\leq (x|y)_b+10\delta.
\ee

Similarly, for $x\in X$ and $\eta\in \partial_\infty X\setminus\{\xi\}$, the Gromov product $(x|\eta)_b$ of $x$ and $\eta$ based at $b$ is defined by
$$(x|\eta)_b= \inf \big\{ \liminf_{i\rightarrow \infty}(x|z_i)_b\;|\; \{z_i\}\in \eta\big\}.$$
For points $\eta,\zeta\in  \partial_\infty X\setminus\{\xi\}$, we define their Gromov product based at $b$ by
$$(\eta|\zeta)_b=\inf\big\{\liminf_{i\to\infty} (x_i|y_i)_b\;|\; \{x_i\}\in\eta , \{y_i\}\in\zeta\big\}.$$

\begin{lem}$($\cite[Lemma 3.2.4]{BuSc}$)$\label{BuSc-Lem3.2.4}
Suppose that $(X, d)$ is a $\delta$-hyperbolic space. For any $\eta$, $\zeta\in \partial_\infty X\setminus\{\xi\}$, and for any Gromov sequences $\{x_i\}\in \eta$ and $\{y_i\}\in \zeta$, we have
$$(\eta|\zeta)_b\leq \liminf_{i\to \infty}(x_i|y_i)_b \leq \limsup_{i\to \infty}(x_i|y_i)_b\leq (\eta|\zeta)_b+44\delta.$$
\end{lem}

For $\varepsilon>0$ with $e^{22\varepsilon\delta}\leq 2$, we define
$$\rho_{b,\varepsilon}(\eta,\zeta)= e^{-\varepsilon(\xi_1|\xi_2)_b}\;\;\;\;\;\;\;\mbox{for all}\;\eta,\zeta\in \partial_\infty X\setminus\{\xi\}.$$
Then for $i=1,2,3$ with $\xi_i\in \partial_\infty X\setminus\{\xi\}$, we have
$$\rho_{b,\varepsilon}(\xi_1,\xi_2)\leq e^{22\varepsilon\delta} \max\{\rho_{b,\varepsilon}(\xi_1,\xi_3),\rho_{b,\varepsilon}(\xi_3,\xi_2)\}.$$
We now define
$$d_{b,\varepsilon}(\eta,\zeta):=\inf\Big\{\sum_{i=1}^{n} \rho_{b,\varepsilon} (\zeta_{i-1},\zeta_i)\;|\;n\geq 1,\zeta=\zeta_0,\zeta_1,\ldots,\zeta_n=\eta\in\partial_\infty X\setminus\{\xi\}\Big\}.$$
By \cite[Lemma $3.3.3$]{BuSc}, it follows that $(\partial_\infty X\setminus\{\xi\}, d_{b,\varepsilon})$ is a metric space such that
\be\label{eq-b}
 \frac{1}{2}\rho_{b,\varepsilon} \leq d_{b,\varepsilon}\leq \rho_{b,\varepsilon}.
 \ee
Then $d_{b,\varepsilon}$ is called a Hamenst\"adt metric on the punctured space $\partial_\infty X\setminus\{\xi\}$ based at $\xi$ with the parameter $\varepsilon$.

%%%%%%%%%%%%%%%%%
%%%%%%%%%%%%%%%%%
\section{Gromov hyperbolic domains}\label{sec-3}
%%%%%%%%%%%%%%%%%
%%%%%%%%%%%%%%%%%

Recall that a domain $G\subsetneq \mathbb{R}^n$ is  Gromov hyperbolic if $(G, k)$ is $\delta$-hyperbolic for some $\delta\geq 0$, where $k$ denotes the quasihyperbolic metric of $G$. In this part, we assume that $G$ is a $\delta$-hyperbolic domain and $\partial_\infty G$ is its Gromov boundary. By \cite[Proposition 2.8]{BHK}, the metric space $(G, k)$ is proper and geodesic.

Fix $\xi_0\in\partial_\infty G$ and a Busemann function $b=b_{\xi_0,o}:(G, k)\to \mathbb{R}$ based at $\xi_0$, where $o\in G$. Consider the family of conformal deformations of $(G, k)$ by the densities
\be\label{b-0} \rho_\varepsilon(x)=e^{-\varepsilon b(x)},\;\;\;\;\;\;\;\;\;\;\;\;\;\;\;\;\;\;\;\;\;\mbox{for all}\;\;\;\varepsilon>0.\ee
We denote the resulting metric spaces by $G_\varepsilon=(G,d_\varepsilon)$. Thus $d_\varepsilon$ is a metric on $G$ defined by
\be\label{b-0b}d_\varepsilon(x,y)=\inf \int_\gamma \rho_\varepsilon \;ds_k,\ee
 where the infimum is taken over all rectifiable curves in $G$ joining the points $x$ and $y$ and $ds_k$ is the arc-length element with respect to the metric $k$.

Denote the metric completion and boundary of $(G_\varepsilon, d_\varepsilon)$ by $\overline{G}_\varepsilon$ and
$$\partial_\varepsilon G:=\partial G_\varepsilon=\overline{G}_\varepsilon\setminus G_\varepsilon,$$
respectively. By (\ref{b-1}), for all $x,y\in G$,
$$|b(x)-b(y)|\leq k(x,y)+10\delta.$$
This guarantees the following {\it Harnack type inequality}:
\be\label{eq-hi} e^{-10\varepsilon\delta}e^{-\varepsilon k(x,y)}\leq \frac{\rho_\varepsilon(x)}{\rho_\varepsilon(y)}\leq e^{10\varepsilon\delta}e^{\varepsilon k(x,y)},\ee
for $x,y\in G$ and $\varepsilon>0$.

Next, let $\mathcal{L}_n$ be the $n$-dimensional Lebesgue measure on $\mathbb{R}^n$. For any Borel set $E\subset G$, we define
\be\label{b-1b}\mu_\varepsilon(E)=\int_E \bigg(\frac{\rho_\varepsilon(x)}{d(x)}\bigg)^n \,d\mathcal{L}_n(x).\ee

For later use, we establish an auxiliary lemma about Gromov hyperbolic domains in $\mathbb{R}^n$.

\begin{lem}\label{main lem-2}
Let $G\subsetneq \mathbb{R}^n$ be a $\delta$-hyperbolic domain, $o\in G$, $\xi_0\in \partial_\infty G$ and $b=b_{\xi_0, o}\in \mathcal{B}(\xi_0)$ a Busemann function on $(G, k)$. Let $(G_\varepsilon, d_\varepsilon, \mu_\varepsilon)$ be the deformed metric measure space which is induced by the density $\rho_\varepsilon$ as in \eqref{b-0}. Then for some $0<\varepsilon=\varepsilon(\delta)$, we have the following:
\begin{enumerate}
\item\label{lem2-1} $(G, k)$ is a proper and geodesic metric space;
\item\label{lem2-2} $(G, k)$ is $K$-roughly starlike with respect to each point of $\partial_\infty G$;
\item\label{lem2-3} $(G_\varepsilon, d_\varepsilon)$ is an unbounded $A$-uniform space, and for each pair of points $x, y\in G$, we have
    $$M^{-1}k(x, y)\leq k_\varepsilon(x, y)\leq M k(x, y),$$
    where $k_\varepsilon$ is the quasihyperbolic metric of $(G_\varepsilon, d_\varepsilon)$;
\item\label{lem2-3a} For all $x, y\in G$,
$$C_\delta^{-1} d_\varepsilon(x, y)\leq \varepsilon^{-1} e^{-\varepsilon(x|y)_b}\min\{1,\varepsilon k(x, y)\}\leq C_\delta d_\varepsilon(x, y);$$
\item\label{lem2-4} $(G_\varepsilon, d_\varepsilon, \mu_\varepsilon)$ is a $n$-Loewner with a control function $\phi=\phi(n, \delta)$, and a locally Ahlfors $n$-regular metric measure space;
\item\label{lem2-3b} Let $\Gamma$ be a family of curves in $G$. Then
$${\rm mod}_n(\Gamma; G)={\rm mod}_n(\Gamma; G_\varepsilon);$$
%\item\label{lem2-5} There is a natural identification $G\cup \partial_\infty G \xrightarrow{\varphi}\overline G_\varepsilon\cup\{\infty\}$;
\item\label{lem2-6} There is a natural identification $\varphi:\, \partial_\infty G\to \partial_\varepsilon G \cup \{\infty\}$ with $\varphi(\xi_0)=\infty$ such that
$$\varphi:\,(\partial_\infty G \setminus\{\xi_0\}, d_{b, \varepsilon_0})\to (\partial_\varepsilon G, d_\varepsilon)$$
is $\eta_0$-quasisymmetric, where $d_{b,\varepsilon_0}$ is a Hamenst\"adt metric based at a Busemann function $b\in \mathcal{B}(\xi_0)$ with parameter $\varepsilon_0$;
\item\label{lem2-7}The identity map $(G, |\cdot|)\to (G_\varepsilon, d_\varepsilon)$ and its inverse map are $q$-locally $\theta_0$-quasisymmetric.
\end{enumerate}
The parameters $K$, $A$, $M$, $C_\delta \geq 1$, $q\in (0, 1)$, $\eta_0$, $\varepsilon_0$, and $\theta_0$ depend only on $\delta$.
\end{lem}

\bpf
(\ref{lem2-1}) See \cite[Proposition 2.8]{BHK}.

(\ref{lem2-2}) See \cite[Lemma 2.3(b)]{ZR}.

(\ref{lem2-3}) See \cite[Theorem 4.1 and Lemma 5.5]{Z20}.

(\ref{lem2-3a}) See \cite[Lemma 5.1]{Z20}.

Next we  give the proofs of Statements (\ref{lem2-4}) $\sim$
%, (\ref{lem2-4}), (\ref{lem2-6}), and
(\ref{lem2-7}) of Lemma \ref{main lem-2}.

\subsection{Proof of Statement (\ref{lem2-4})}

Consider the metric measure space $(G, \ell_G, \mathcal{L}_n)$, where $\ell_G$ is the length metric of $G$ associated with the Euclidean metric $|\cdot|$. Note that $(G, \ell_G, \mathcal{L}_n)$ is a locally compact, non-complete, length metric measure space. It is known that every Euclidean ball is $n$-Loewner with some control function $\phi_n$ depending only on $n$. Because $B(x, r)=B_{\ell_G}(x, r)=\{y\in G\,|\,\ell_G(y, x)<r\}$ for all $0<r<d(x)$, we know that $(G, \ell_G, \mathcal{L}_n)$ is also locally $n$-Loewner.

Let $0<\varepsilon=\varepsilon(\delta)$.
Now consider the conformal density
\be\label{eq4-1}
\rho(x)=\frac{\rho_\varepsilon(x)}{d(x)}=\frac{e^{-\varepsilon b(x)}}{d(x)},
\ee
for $x\in G$. Let $(G, d_\rho, \mu_\rho)$ be the conformal deformation of $(G, \ell_G, \mathcal{L}_n)$ induced by the density $\rho$. That is,
$$d_\rho(x, y)=\inf \int_\gamma \rho(x)\, ds,$$
 and
 $$\mu_\rho(E)=\int_E \rho(x)^n \,d\mathcal{L}_n(x),$$
 where the infimum is taken over all rectifiable curves $\gamma$ in $G$ with end points $x$ and $y$, $ds$ is the arc-length element with respect to the metric $|\cdot|$, and $E$ is a Borel set. Note that
 $$(G, d_\rho, \mu_\rho)=(G, d_\varepsilon, \mu_\varepsilon).$$

For $z\in G$ and all points $x, y\,\in B(z, d(z)/2)$, we see from Lemma \ref{Vai-4 thm3.9} that
$$\frac{1}{3}d(y)\leq d(x)\leq 3d(y)\,\,\,\mbox{and}\,\,\, k(x, y)\leq 2.$$
This together with the Harnack inequality \eqref{eq-hi}, yields
\be\label{eq4-2}
\frac{1}{C_1}\leq\frac{\rho(x)}{\rho(y)}=\frac{\rho_\varepsilon(x)}{\rho_\varepsilon(y)}
\frac{d(y)}{d(x)}\leq C_1,
\ee
where $C_1=C_1(\delta)$.

By Statement (\ref{lem2-2}), $(G, k)$ is $K$-roughly starlike with respect to $\xi_0\in \partial_\infty G$, where $K$ depends only on $\delta$. Hence it follows from \cite[Lemma 5.3]{Z20} that for all $x\in G$,
$$\frac{1}{C_2}\rho_\varepsilon(x) \leq d_\varepsilon(x)\leq C_2 \rho_\varepsilon(x),$$
where $C_2=C_2(\delta)$, which implies
\be\label{eq4-4}
\frac{1}{C_2}\leq \frac{d_\rho(x)}{d(x)\rho(x)}=\frac{d_\varepsilon(x)}{\rho_\varepsilon(x)}\leq C_2.
\ee
 Combining this with \cite[Theorem 6.39]{BHK}, we observe that $(G, d_\varepsilon, \mu_\varepsilon)$ is a locally $n$-Loewner space with the control function depending only on $n$ and $\delta$. Thus from \cite[Theorem 6.4]{BHK}, we know that $(G, d_\varepsilon, \mu_\varepsilon)$ is actually $n$-Loewner with the control function $\phi=\phi(n, \delta)$, because $(G, d_\varepsilon)$ is $A$-uniform with $A=A(\delta)$ by using Statement (\ref{lem2-3}).

It remains to prove that $(G, d_\varepsilon, \mu_\varepsilon)$ is locally Ahlfors $n$-regular. Take a number $0<\lambda<1/(32A^2 M)$, where $A$ and $M$ are the constants of Statement (\ref{lem2-3}). Fix $x\in G$ and $0<r<\lambda d_\varepsilon(x)$.
We show that
\be\label{eq4-3}
B\Big(x, \frac{d(x)r}{4e M d_\varepsilon(x)}\Big)\subset B_\varepsilon(x, r)=\{u\in G\,|\,d_\varepsilon(u, x)<r\}\subset B\Big(x,\frac{16 A^2 M d(x)r}{d_\varepsilon(x)}\Big).
\ee

First, we prove the first inclusion. For all $y\in B\big(x, d(x)r/(4 e M d_\varepsilon (x))\big)$, we see from the choices of $\lambda$ and $r$ that
\beq\label{eq4-3b}
|x-y|< \frac{d(x)r}{4e M d_\varepsilon(x)}\leq \frac{d(x)\lambda}{4e M}< \frac{ d(x)}{2M}.
\eeq
Then by Lemma \ref{Vai-4 thm3.9}, we have
\be\label{eq4-3c}
k(x, y)\leq 2\frac{|x-y|}{d(x)}.
\ee
It follows from \eqref{eq4-3b}, \eqref{eq4-3c}, Statement (\ref{lem2-3}), and Lemma \ref{BHK 2.4} that
\beq\label{eq4-3g}
\log\bigg(1+\frac{d_\varepsilon(x, y)}{d_\varepsilon(x)}\bigg)\leq k_\varepsilon(x, y)\leq M k(x, y)\leq 2M\frac{|x-y|}{d(x)}\leq 1,
\eeq
which shows that
$$\frac{1}{e}\frac{d_\varepsilon(x, y)}{d_\varepsilon(x)}\leq \log\bigg(1+\frac{d_\varepsilon(x, y)}{d_\varepsilon(x)}\bigg),$$
because $\log(1+t)\geq (1/e) t$ for $0\leq t\leq e-1$. Hence, we obtain from \eqref{eq4-3g} that
$$\frac{d_\varepsilon(x, y)}{d_\varepsilon(x)}\leq 2eM\frac{|x-y|}{d(x)}.$$
Furthermore, because $y\in B\big(x, d(x)r/(4 e M d_\varepsilon (x))\big)$, we have
$$d_\varepsilon(x, y)\leq \frac{r}{2}<r.$$

Next, we prove the second inclusion in \eqref{eq4-3}. For all $y\in B_\varepsilon(x, r)$,
we have
\be\label{eq4-3e}
d_\varepsilon(y)\geq d_\varepsilon(x)-d_\varepsilon(x, y)\geq (1-\lambda)d_\varepsilon(x)\geq \frac{1}{2}d_\varepsilon(x).
\ee
Then we see from  Theorem \ref{BHK 2.4} and Statement (\ref{lem2-3}) that
\beq\label{eq4-3d}
k(x, y)&\leq& M k_\varepsilon(x, y)\nonumber\\
&\leq& 4A^2M\log\bigg(1+\frac{d_\varepsilon(x, y)}{\min\{d_\varepsilon(x), d_\varepsilon(y)\}}\bigg)\nonumber\\
&\leq& 8A^2 M \frac{d_\varepsilon(x, y)}{d_\varepsilon(x)}.
\eeq
Because $y\in B_\varepsilon(x, r)$ and $r<\lambda d_\varepsilon(x)$, \eqref{eq4-3d} implies that $k(x, y)\leq 1.$
Hence, it follows from Lemma \ref{Vai-4 thm3.9} and \eqref{eq4-3d} that
$$\frac{|x-y|}{d(x)}\leq 2k(x, y)\leq 16A^2 M \frac{d_\varepsilon(x, y)}{d_\varepsilon(x)},$$
and, therefore,
$$|x-y|\leq \frac{16A^2 M d_\varepsilon(x, y)d(x)}{d_\varepsilon(x)}<\frac{16A^2 M d(x)r}{d_\varepsilon(x)}.$$
Hence \eqref{eq4-3} holds true.

Because $r<\lambda d_\varepsilon(x)$ and $\lambda<1/(32A^2 M)$, we know that
\be\label{eq4-3f}
B\Big(x, \frac{16 A^2 M d(x)r}{d_\varepsilon(x)}\Big)\subset B\Big(x, \frac{d(x)}{2}\Big)\subset G.
\ee
Recall that
$$\mu_\varepsilon \big(B_\varepsilon(x, r)\big)=\mu_\rho\big(B_\varepsilon(x, r)\big)=\int_{B_\varepsilon(x, r)} \rho(y)^n\,d\mathcal{L}_n(y),$$
Therefore, we obtain from \eqref{eq4-2} and \eqref{eq4-3f} that
\beq\label{eq4-3a}
\frac{1}{C_1^n} \rho(x)^n \mathcal{L}_n \big(B_\varepsilon(x, r)\big)\leq \mu_\varepsilon\big(B_\varepsilon(x, r)\big)\leq C_1^n \rho(x)^n \mathcal{L}_n\big(B_\varepsilon(x, r)\big).
\eeq
Because $(G, |\cdot|, \mathcal{L}_n)$ is locally Ahlfors $n$-regular with the constant $C_0$, we obtain
\beqq
\mu_\varepsilon\big(B_\varepsilon(x, r)\big)&\leq& C_1^n\rho(x)^n\mathcal{L}_n\bigg(B\Big(x, \frac{16 A^2 M d(x) r}{d_\varepsilon (x)}\Big)\bigg)\,\,\,(\mbox{By}\,\,\,\,\, \eqref{eq4-3}\,\,\,\mbox{and}\,\,\,\eqref{eq4-3a})\\
&\leq& C_0C_1^n\rho(x)^n\frac{(16 A^2 M d(x))^n}{d_\varepsilon(x)^n}r^n\\
&\leq& C_0(16 C_1 C_2 A^2 M)^n r^n \,\,\quad\quad\quad\quad\quad\quad\,\,\,(\mbox{By}\,\,\, \eqref{eq4-4})\\
&=& C_3 r^n.
\eeqq
On the other hand, we have
\beqq
\mu_\varepsilon\big(B_\varepsilon(x, r)\big)&\geq& \frac{1}{C_1^n}\rho(x)^n \mathcal{L}_n\bigg(B\Big(x, \frac{d(x)r}{4eMd_\varepsilon(x)}\Big)\bigg)\,\,\,\,\,\,\,\,(\mbox{By}\,\,\,\,\, \eqref{eq4-3}\,\,\,\mbox{and}\,\,\,\eqref{eq4-3a})\\
&\geq& \frac{1}{C_0C_1^n}\frac{d(x)^n}{(4 e Md_\varepsilon(x))^n}\rho(x)^n r^n\\
&\geq& \frac{1}{C_0(4 e C_1 C_2 M)^n}r^n \,\,\quad\quad\quad\quad\quad\quad\,\,\,\,\,\,\,\,(\mbox{By}\,\,\, \eqref{eq4-4})\\
&\geq& \frac{1}{C_3}r^n,
\eeqq
as desired.
\qed

\subsection{Proof of Statement (\ref{lem2-3b})} Let $\Gamma$ be a family of curves in $G$. We need to show that ${\rm mod}_n(\Gamma; G)={\rm mod}_n(\Gamma; G_\varepsilon)$.

Let $\rho$ be the density function defined as in \eqref{eq4-1} and $\rho_\varepsilon$ be as in \eqref{b-0}. It follows from \eqref{b-0b} that
$$ds_\varepsilon=\rho_\varepsilon(x)\,ds_k=\frac{\rho_\varepsilon(x)}{d(x)}\,ds=\rho(x)\,ds,$$
where $ds_\varepsilon$ is the arc-length element in $(G, d_\varepsilon, \mu_\varepsilon)$.
By \eqref{b-1b}, we have
  $$d\mu_\varepsilon(x)=\rho(x)^n\,d\mathcal{L}_n(x),$$
 where  $\mu_\varepsilon$ is the measure element in $(G, d_\varepsilon, \mu_\varepsilon)$.

Now, for any curve $\gamma\in \Gamma$, and for all Borel functions $\rho_1: (G, |\cdot|, \mathcal{L}_n)\to [0, +\infty]$ satisfying $\int_{\gamma} \rho_1(x)\,ds\geq 1$, we have
$$\int_{\gamma} \rho_2(x)\,ds_{\varepsilon}=\int_{\gamma}\rho_2(x)\rho(x)\,ds=\int_{\gamma} \rho_1(x)\,ds\geq 1,$$
where $\rho_2(x)=\rho_1(x)/\rho(x)$ and $\rho_2(x):(G_\varepsilon, d_\varepsilon, \mu_\varepsilon)\to [0, +\infty]$ is a Borel function. Thus we see from the definition of modulus that
$${\rm mod}_n(\Gamma; G_\varepsilon)\leq \int_{G}\rho_2(x)^n\,d\mu_\varepsilon(x)=\int_G \big(\rho_2(x)\rho(x)\big)^n\,d\mathcal{L}_n(x)=\int_G \rho_1(x)^n\,d\mathcal{L}_n(x).$$
Hence, by the arbitrariness of $\rho_1$, we obtain
$${\rm mod}_n(\Gamma; G_\varepsilon)\leq {\rm mod}_n(\Gamma; G).$$
Similarly, we get
$${\rm mod}_n(\Gamma; G)\leq {\rm mod}_n(\Gamma; G_\varepsilon),$$
as desired.
\qed

\subsection{Proof of Statement (\ref{lem2-6})}
By \cite[Lemma 5.2]{Z20}, we see that there is a well defined identification
$$\varphi:\, \partial_\infty G\to \partial_\varepsilon G\cup \{\infty\}$$
with $\varphi(\xi_0)=\infty$ for some $\xi_0\in \partial_\infty G$. Let $b=b_{\xi_0, o}\in \mathcal{B}(\xi_0)$ be a Busemann function and $d_{b, \varepsilon_0}$ a Hamenst\"adt metric based at $b$ with the parameter $\varepsilon_0=\varepsilon_0(\delta)$. If $f$ is $\theta$-quasisymmetric, then $f^{-1}$ is $\theta_1$-quasisymmetric \cite[Theorem 2.2]{TV}, where $\theta_1(t)=1/\theta^{-1}(1/t)$. Thus it suffices to prove that the map $\varphi^{-1}$ is $\eta_0$-quasisymmetric.

For all distinct points $\xi, \eta, \zeta\in \partial_\infty G\setminus\{\xi_0\}$, choose Gromov sequences $\{x_n\}\in \xi$, $\{y_n\}\in \eta$, and $\{z_n\}\in \zeta$, respectivly. By the definition of $\varphi$, we know that $\{x_n\}\xrightarrow{d_\varepsilon}\varphi(\xi)$, $\{y_n\}\xrightarrow{d_\varepsilon}\varphi(\eta)$, and $\{z_n\}\xrightarrow{d_\varepsilon}\varphi(\zeta)$, as $n\to \infty$, respectively. We see from Lemma \ref{z0} that
$$\liminf_{n\to \infty}(x_n|y_n)_o=\liminf_{n\to \infty} \frac{1}{2}\big(k(x_n,o)+k(y_n,o)-k(x_n, y_n)\big)\leq (\xi|\eta)_o+2\delta.$$
Because $(x_n|x_n)_o=k(x_n, o)\to +\infty$ as $n\to \infty$,
this shows that
$$\liminf\limits_{n\to \infty}k(x_n, y_n)=+\infty.$$
 Then without loss of generality, we may assume that for all $n$,
\be\label{eq4-5}
\varepsilon k(x_n, y_n)\geq 1.
\ee
It follows from Statement (\ref{lem2-3a}) and \eqref{eq4-5}  that
\beq\label{eq4-5a}
\frac{d_{\varepsilon}(x_n, y_n)}{d_{\varepsilon}(x_n, z_n)}&\geq& \frac{1}{C_\delta^2} \frac{e^{-\varepsilon (x_n|y_n)_b}\min\{1, \varepsilon k(x_n, y_n)\}}{e^{-\varepsilon (x_n|z_n)_b}\min\{1, \varepsilon k(x_n, z_n)\}}\nonumber\\
&\geq& \frac{1}{C_\delta^2}e^{\varepsilon( (x_n|z_n)_b- (x_n|y_n)_b)},
\eeq
where $C_\delta$ is the constant of Statement (\ref{lem2-3a}). Now by Lemma \ref{BuSc-Lem3.2.4}, it follows that
\beqq
(\xi|\eta)_b-(\xi|\zeta)_b\leq \liminf\limits_{n\to \infty}\big((x_n|z_n)_b-(x_n|y_n)_b\big)+44\delta.
\eeqq
Hence, by \eqref{eq-b} and \eqref{eq4-5a}, we have
\beqq
\frac{d_{b, \varepsilon_0}(\xi, \eta)}{d_{b, \varepsilon_0}(\xi, \zeta)}&\leq& 2\frac{\rho_{b, \varepsilon_0}(\xi, \eta)}{\rho_{b, \varepsilon_0}(\xi, \zeta)} =2e^{\varepsilon_0(\xi|\zeta)_b-\varepsilon_0(\xi|\eta)_b}\\
&\leq&\liminf\limits_{n\to \infty}2e^{\varepsilon_0(x_n|z_n)_b-\varepsilon_0(x_n|y_n)_b+44\delta \varepsilon_0}\\
&\leq&2e^{44\delta\varepsilon_0}C_\delta^{\frac{2\varepsilon_0}{\varepsilon}}\liminf\limits_{n\to\infty}
\bigg(\frac{d_{\varepsilon}(x_n, y_n)}{d_{\varepsilon}(x_n, z_n)}\bigg)^{\frac{\varepsilon_0}{\varepsilon}}\\
&\leq& C(\delta)\bigg(\frac{d_\varepsilon(\varphi(\xi), \varphi(\eta))}{d_\varepsilon(\varphi(\xi), \varphi(\zeta))}\bigg)^{\frac{\varepsilon_0}{\varepsilon}}.
\eeqq
\qed

\subsection{Proof of Statement (\ref{lem2-7})}
%%%%%%%%%%%%%%%%%%%
%%%%%%%%%%%%%%%%%%%
First, we prove that the identity map $(G, |\cdot|)\to (G, d_\varepsilon)$ is $q_1$-locally $\theta_0$-quasisymmetric with $q_1=(\log\frac{3}{2})/8M$ and $\theta_0(t)=64 A^2 M^2 t$, where $A$ and $M$ are the constants of Statement (\ref{lem2-3}). For all $x_0\in G$ and all $x, y, z\in B(x_0, q_1d(x_0))$, we have
$$\max \{|x-y|,|x-z|\}\leq 2q_1d(x_0)\leq \frac{2q_1}{1-q_1}d(x)\leq 4q_1d(x).$$
By Lemma \ref{Vai-4 thm3.9}, it follows that
\be\label{eq4-6}
\frac{1}{2}\frac{|x-y|}{d(x)}\leq k(x, y)\leq 2\frac{|x-y|}{d(x)}\leq 8q_1,
\ee
and
\be\label{eq4-7}
\frac{1}{2}\frac{|x-z|}{d(x)}\leq k(x, z)\leq 2\frac{|x-z|}{d(x)}\leq 8q_1.
\ee
Thus, by Theorem \ref{BHK 2.4} and Statement (\ref{lem2-3}), we have
$$\log\bigg(1+\frac{d_\varepsilon(x, y)}{\min\{d_\varepsilon(x), d_\varepsilon(y)\}}\bigg)\leq k_\varepsilon(x, y)\leq Mk(x, y)\leq 8M q_1=\log \frac{3}{2},$$
and therefore,
$$d_\varepsilon(x, y)\leq \frac{1}{2}d_\varepsilon(x).$$
This ensures that
\be\label{eq4-7a}
k_\varepsilon(x, y)\geq \log\bigg(1+\frac{d_\varepsilon(x, y)}{d_\varepsilon(x)}\bigg)\geq \frac{d_\varepsilon(x, y)}{2d_\varepsilon(x)},
\ee
because $\log(1+t)\geq t/2$ for all $0\leq t\leq 1$.
Similarly, we compute that $d_\varepsilon(x, z)\leq \frac{1}{2}d_\varepsilon(x)$ and $d_\varepsilon(z)\geq \frac{1}{2}d_\varepsilon(x)$. By \eqref{eq4-6} and \eqref{eq4-7a}, we obtain
\be\label{eq4-7b}
d_\varepsilon(x, y)\leq 2k_\varepsilon(x, y)d_\varepsilon(x)\leq 2 M k(x, y)d_\varepsilon(x)\leq 4 M \frac{d_\varepsilon(x)}{d(x)}|x-y|.
\ee
Furthermore, by Statement (\ref{lem2-3}), we know that $(G, d_\varepsilon)$ is $A$-uniform, and we infer from Theorem \ref{BHK 2.4} and the fact $d_\varepsilon(z)\geq \frac{1}{2}d_\varepsilon(x)$ that
$$k_\varepsilon(x, z)\leq 4A^2 \log\bigg(1+\frac{d_\varepsilon(x, z)}{\min\{d_\varepsilon(x), d_\varepsilon(z)\}}\bigg)\leq 8A^2\frac{d_\varepsilon(x,z)}{d_\varepsilon(x)}.$$
By \eqref{eq4-7}, we observe that
\be\label{eq4-7c}
d_\varepsilon(x, z)\geq \frac{1}{8A^2}k_\varepsilon(x, z)d_\varepsilon(x)\geq \frac{1}{8A^2 M} k(x, z)d_\varepsilon(x)\geq \frac{1}{16 A^2 M}\frac{d_\varepsilon(x)}{d(x)}|x-z|.
\ee
Therefore, we obtain from \eqref{eq4-7b} and \eqref{eq4-7c} that
$$\frac{d_\varepsilon(x, y)}{d_\varepsilon(x, z)}\leq 64 A^2 M^2\frac{|x-y|}{|x-z|}.$$

It remains to show that the inverse map of the identity $(G, |\cdot|)\to (G, d_\varepsilon)$ is $q_2$-locally $\theta_0$-quasisymmetric with $q_2=1/(32A^2 M)$. For all $x_0\in G$ and all $x, y, z\in B_\varepsilon(x_0, q_2d_\varepsilon(x_0))$, by the triangle inequality, we get
\be\label{eq4-8}
\max\{d_\varepsilon(x, y), d_\varepsilon(x, z)\}\leq 2q_2d_\varepsilon(x_0)\leq \frac{2q_2}{1-q_2}d_\varepsilon(x)\leq 4q_2d_\varepsilon(x),
\ee
which implies
\be\label{eq4-8a}\min\{d_\varepsilon(y), d_\varepsilon(z)\}\geq \frac{1}{2}d_\varepsilon(x).
\ee
Because $(G, d_\varepsilon)$ is $A$-uniform, it follows from \eqref{eq4-8a} and  Theorem \ref{BHK 2.4} that
\be\label{eq4-9}
k_\varepsilon(x, y)\leq 4A^2\log\bigg(1+\frac{d_\varepsilon(x, y)}{\min\{d_\varepsilon(x), d_\varepsilon(y)\}}\bigg)\leq 8A^2 \frac{d_\varepsilon(x, y)}{d_\varepsilon(x)}\leq 32 A^2 q_2.
\ee
Similarly, we see that
\be\label{eq4-10}
k_\varepsilon(x, z)\leq  32A^2 q_2.
\ee
On the other hand, Theorem \ref{BHK 2.4} gives us
\be\label{eq4-10a}
k_\varepsilon(x, z)\geq \log\bigg(1+\frac{d_\varepsilon(x, z)}{d_\varepsilon(x)}\bigg)\geq \frac{d_\varepsilon(x, z)}{2d_\varepsilon(x)},
\ee
here we have used the fact that $\log(1+t)\geq t/2$ for all $0\leq t\leq 1$.
Therefore, we see from Statement (\ref{lem2-3}) that
\beqq
\max\{k(x, y), k(x, z)\}\leq M \max\{k_\varepsilon(x, y), k_\varepsilon(x, z)\}\leq 32A^2 M q_2= 1.
\eeqq
Now, by  \eqref{eq4-8}, Lemma \ref{Vai-4 thm3.9}, and Statement (\ref{lem2-3}), it follows that
%\be\label{eq4-9}
%\frac{1}{2}\frac{d_\varepsilon(x, y)}{d_\varepsilon(x)}\leq \frac{1}{2}\frac{d_\varepsilon(x, y)}{\min\{d_\varepsilon(x), d_\varepsilon(y)\}}\leq k_\varepsilon(x, y)\leq 8A^2 \frac{d_\varepsilon(x, y)}{d_\varepsilon(x)}\leq 16 A^2 q_2,
%\ee
%and
%\be\label{eq4-10}
%\frac{1}{2}\frac{d_\varepsilon(x, z)}{d_\varepsilon(x)}\leq \frac{1}{2}\frac{d_\varepsilon(x, z)}{\min\{d_\varepsilon(x), d_\varepsilon(z)\}}\leq k_\varepsilon(x, z)\leq 8A^2 \frac{d_\varepsilon(x, z)}{d_\varepsilon(x)}\leq 16 A^2 q_2,
%\ee
%because $t/2\leq \log(1+t)\leq t$ for $0\leq t\leq 1$. Moreover, we see from Statement (\ref{lem2-3}) that
%$$\max\{k(x, y), k(x, z)\}\leq M \max\{k_\varepsilon(x, y), k_\varepsilon(x, z)\}\leq 16A^2 M q_2= 1.$$
$$|x-y|\leq 2d(x)k(x, y)\leq 2 M d(x) k_\varepsilon(x, y)\leq \frac{16 A^2 M d(x)}{d_\varepsilon(x)}d_\varepsilon(x, y).$$
Similarly, using \eqref{eq4-10a}, we have
$$|x-z|\geq \frac{1}{2}d(x)k(x, z)\geq \frac{1}{2 M}k_\varepsilon(x, z)d(x)\geq \frac{d(x)}{4 M d_\varepsilon(x)}d_\varepsilon(x, z).$$
Hence,
$$\frac{|x-y|}{|x-z|}\leq 64 A^2 M^2\frac{d_\varepsilon(x, y)}{d_\varepsilon(x, z)}.$$
\epf

The proof of Lemma \ref{main lem-2} is complete.
\qed
%%%%%%%%%%%%%%%%%%%
%%%%%%%%%%%%%%%%%%%
\section{Proofs of Theorems \ref{main thm-1} and \ref{thm-4}}\label{sec-4}
%%%%%%%%%%%%%%%%%%%
%%%%%%%%%%%%%%%%%%%
In this section, we study geometric properties of Gromov hyperbolic domains. The main objective is to show Theorem \ref{main thm-1} and Theorem \ref{thm-4}. We fist recall the following definition:
\bdefe
For $c\geq 1$, a metric space $X$ is called {\it $C_{bt}$-bounded turning} if for each pair of points $x, y\in X$ can be joined by a curve $\gamma$ with $\diam(\gamma)\leq C_{bt}d(x, y)$.
\edefe

\subsection*{Proof of Theorem \ref{main thm-1}}
Assume that $G\subsetneq \mathbb{R}^n$ is a $\delta$-hyperbolic domain. Let $(G, d_\varepsilon, \mu_\varepsilon)$ be the metric measure space induced by the density \eqref{b-0} on $(G, |\cdot|, \mathcal{L}_n)$, for some constant $0<\varepsilon=\varepsilon(\delta)$. Note that the second statement follows immediately from the first statement and the definition of the bounded turning condition. So we only need to prove the first assertion.

Fix $x, y\in G$. By Lemma \ref{main lem-2}(\ref{lem2-1}), there is a quasihyperbolic geodesic $\alpha=[x, y]_k$ in $G$ between $x$ and $y$. We observe from the proof of \cite[Theorem 4.1]{Z20} that $\alpha$ is an $A_1$-uniform curve in the deformed space $(G, d_\varepsilon)$ with $A_1=A_1(\delta)$. Thus it follows that for all $z\in \alpha$,
\be\label{eq5-1}
\min\{\diam_\varepsilon(\alpha[x, z]), \diam_\varepsilon(\alpha[z, y])\}\leq A_1d_\varepsilon(z),
\ee
where $\diam_\varepsilon$ means the diameter in the metric $d_\varepsilon$ and $d_\varepsilon(x)=\dist_\varepsilon(x, \partial_\varepsilon G)$ for all $x\in G$.
Without loss of generality, we may assume that $\diam_\varepsilon(\alpha[x, z])\leq \diam_\varepsilon(\alpha[z, y])$. Hence, by \eqref{eq5-1}, we have
\be\label{eq5-2}
\diam_\varepsilon(\alpha[x, z])\leq A_1d_\varepsilon(z).
\ee
%It follows from Lemma \ref{main lem-2}(\ref{lem2-7}) that the identification $(G, |\cdot|)\to G_\varepsilon$ and its inverse map are both $q$-locally $\theta_0$-quasisymmetric.
Let $\gamma$ be a fixed curve in $G$ between $x$ and $y$, and denote
$$R=\diam(\gamma)\,\,\,\mbox{and}\,\,\,r=\frac{q}{A_1}\diam_\varepsilon(\alpha[x, z])\leq qd_\varepsilon(z),$$
where $q\in (0, 1)$ is a constant depending only on $\delta$ such that Statement \eqref{lem2-7} of Lemma \ref{main lem-2} holds.

\begin{figure}[htbp]
\begin{center}
\includegraphics{fig5.1}
\end{center}
\caption{} \label{f1}
\end{figure}

 We observe that there exists a subcurve $\alpha_1$ of $\alpha[x, z]$ that intersects both $S_\varepsilon(z, r/2)$ and $S_\varepsilon(z, r/4)$ and lies in $\overline{B}_\varepsilon(z, r/2)\setminus B_\varepsilon (z, r/4)$, where $S_\varepsilon$, $B_\varepsilon$ and $\overline{B}_\varepsilon$ denote the sphere, open ball, and closed ball in the metric $d_\varepsilon$, respectively. Similarly, there is  a subcurve $\alpha_2$ of $\alpha[x, z]$ which intersects both $S_\varepsilon(z, r/8)$ and $S_\varepsilon(z, r/16)$ and is contained in $\overline{B}_\varepsilon(z, r/8)\setminus B_\varepsilon(z, r/16)$. See Figure \ref{f1}.

 Because $\alpha$ is $A_1$-uniform, by \eqref{eq5-2} and by the choice of $\alpha_1$, we obtain
 $$\diam_\varepsilon(\gamma)\geq d_\varepsilon(x, y)\geq \frac{1}{A_1}\diam_\varepsilon(\alpha[x, y])\geq r\geq \diam_\varepsilon(\alpha_1)\geq \frac{r}{4},$$
 and
 $$\dist_\varepsilon(\alpha_1, \gamma)\leq \diam_\varepsilon(\alpha[x,z])=\frac{A_1 r}{q}.$$
 The above estimates ensure that
 $$\Delta_\varepsilon(\gamma, \alpha_1)=\frac{\dist_\varepsilon (\gamma, \alpha_1)}{\min\{\diam_\varepsilon(\gamma), \diam_\varepsilon(\alpha_1)\}}=\frac{\dist_\varepsilon (\gamma, \alpha_1)}{\diam_\varepsilon(\alpha_1)}\leq \frac{4A_1}{q}:=A_2.$$
Moreover, Lemma \ref{main lem-2}(\ref{lem2-4}) yields that $(G_\varepsilon, d_\varepsilon, \mu_\varepsilon)$ is $n$-Loewner with the control function $\phi=\phi(n, \delta)$. Therefore, we have
\be\label{eq5-3}
\mbox{mod}_n(\gamma, \alpha_1; G_\varepsilon)\geq \phi(\Delta_\varepsilon(\gamma, \alpha_1))\geq \phi(A_2).
\ee
Now, let
$$\dist(x, \alpha_1)=C_1R=C_1\diam(\gamma).$$
We shall find an upper bound for the constant $C_1$. Without loss of generality, we may assume that $C_1\geq 4$. Note that
$$\alpha_1\subset G\setminus B(x,C_1R).$$
Because $\gamma\subset \overline{B}(x, R)$, for any curve $\sigma$  in $G$ joining $\alpha_1$ and $\gamma$, there is a subcurve $\sigma_1\subset \sigma$  connecting $S(x, R)$ and $S(x, C_1 R)$. By the standard modulus estimate \cite[7.5]{Vai75}, we obtain
$$\mbox{mod}_n(\alpha_1, \gamma; G)\leq \mbox{mod}_n(S(x, R), S(x, C_1R); G)\leq \omega_{n-1} (\log C_1)^{1-n},$$
where $\omega_{n-1}$ is $(n-1)$-dimensional surface area of the unit sphere $\mathbb{S}^{n-1}$. On the other hand, it follows from Lemma \ref{main lem-2}(\ref{lem2-3b}) that
$$\mbox{mod}_n(\gamma, \alpha_1; G_\varepsilon)=\mbox{mod}_n(\gamma, \alpha_1; G).$$
This, together with \eqref{eq5-3}, shows that
\beqq
\phi(A_2)&\leq& \mbox{mod}_n(\gamma, \alpha_1; G_\varepsilon)=\mbox{mod}_n(\gamma, \alpha_1; G)\leq \omega_{n-1} (\log C_1)^{1-n},
\eeqq
which yields
$$C_1\leq C(\omega_{n-1}, A_2, \phi)=C(n, \delta).$$
Then there is a point $x_1\in \alpha_1$ such that
\be\label{eq5-4a}
|x_1-x|\leq CR.
 \ee

 A similar computation shows that $\dist(x,\alpha_2)\leq C R$ and there is a point $x_2\in \alpha_2$ satisfying $|x-x_2|\leq C R$. Hence we have
$$|x_1-x_2|\leq |x-x_1|+|x-x_2|\leq 2CR.$$
By the choices of $\alpha_1$ and $\alpha_2$, it follows that $\alpha_1\cup \alpha_2\subset \overline{B}_\varepsilon(z, qd_\varepsilon(z))$, and
\be\label{eq5-4b}
d_\varepsilon(x_1, z)\leq r/2\leq 4d_\varepsilon(x_1, x_2).
\ee
By Lemma \ref{main lem-2}(\ref{lem2-7}), the inverse of the identity map $(G, |\cdot|)\to (G, d_\varepsilon)$ is $q$-locally $\theta$-quasisymmetric with $\theta$ depending only on $\delta$. Then we derive from \eqref{eq5-4b} that
$$|x_1-z|\leq \theta(4)|x_1-x_2|\leq 2\theta(4)C R.$$
Now, by \eqref{eq5-4a}, we obtain that
$$|x-z|\leq |x-x_1|+|x_1-z|\leq \big(1+2\theta(4)\big)C R.$$
Because this holds for all $z\in \alpha$, it guarantees that
$$\diam(\alpha)\leq 2C\big(1+2\theta(4)\big)\diam(\gamma).$$
\qed

Next, we show Theorem \ref{thm-4} by using Theorem \ref{main thm-1}.

%%%%%%%%%%%%%%%%%%%
%%%%%%%%%%%%%%%%%%%
\subsection*{Proof of Theorem \ref{thm-4}}\label{sec-6}
%%%%%%%%%%%%%%%%%%%
%%%%%%%%%%%%%%%%%%%
Because $G$ is a $\delta$-hyperbolic domain, by Theorems \ref{BB03-thm0.1} and \ref{main thm-1}, there are positive constants $C_{gh}$, $C_{sp}$, and $C_{po}$ depending only on $\delta$ and $n$ such that $G$ satisfies the Gehring-Hayman property with the constant $C_{gh}$, the separation property with the constant $C_{sp}$, and the Pommerenke property with the constant $C_{po}$, respectively.

Let $\Sigma$ be a cross-section in $G$ with $G_1$ and $G_2$ the components of $G\setminus \Sigma$, and let $L$ be a quasihyperbolic geodesic in $G$ whose end points lie in $G_1$, and let $x\in L\cap G_2$ be given. Then there is a subarc $L[y, z]$ of $L$ which contains $x$ and is contained in $G_2$ except for its end points $y$ and $z$. Note that $L[y, z]$ is also a quasihyperbolic geodesic of $G$. Choose a curve $\gamma$ in $\Sigma$ connecting $y$ and $z$. See Figure \ref{f2}.

\begin{figure}[htbp]
\begin{center}
\includegraphics{fig1.1}
\end{center}
\caption{} \label{f2}
\end{figure}

We consider two possibilities. First, suppose $B(x, d(x))$ contains one of $\{y, z\}$, say $y$. In this case, denote $\alpha=L[x, y]$. Obviously, $\alpha$ is an arc connecting $x$ to the cross-section $\Sigma$. By using the Gehring-Hayman property, we obtain
\be\label{eq4-9a}
\ell(\alpha)=\ell(L[x, y])\leq C_{gh}\ell_G(x, y)=C_{gh}|x-y|\leq C_{gh}d(x),
\ee
where $\ell_G$ is the length metric of $G$ induced by the Euclidean metric. Furthermore, it follows from \eqref{eq4-9a} and the Pommerenke property that
$$\ell(\alpha)\leq C_{gh}|x-y|\leq C_{gh}\diam(L[y, z])\leq C_{gh}C_{po}\diam(\gamma)\leq C_{gh}C_{po}\diam(\Sigma).$$
Hence, we obtain
$$\ell(\alpha)\leq C_{gh}C_{po}\min\{d(x), \diam(\Sigma)\}.$$

In the remaining possibility, suppose that $y, z\not\in B(x, d(x))$. On the one hand, we see from the separation property that there is $x_0\in \gamma$ such that
\be\label{eq4-9b}
\ell_G(x, x_0)\leq C_{sp}d(x).
\ee
In this case, we let $\alpha$ be a quasihyperbolic geodesic of $G$ which joins $x$ to $x_0$. Now, applying the Gehring-Hayman property, it follows from \eqref{eq4-9b} that
$$\ell(\alpha)\leq C_{gh}\ell_G(x, x_0)\leq C_{gh}C_{sp}d(x).$$
On the other hand, because $y$ and  $z$ are not in $B(x, d(x))$, by the Pommerenke property, we have
$$d(x)\leq |x-y|\leq \diam(L[y, z])\leq C_{po}\diam(\gamma)\leq C_{po}\diam(\Sigma),$$
and, therefore,
$$\ell(\alpha)\leq C_{gh}C_{sp}C_{po} \min\{d(x), \diam(\Sigma)\}.$$
The theorem follows by taking $A=C_{gh}C_{sp}C_{po}$.
\qed

%%%%%%%%%%%%%%%%%%%
%%%%%%%%%%%%%%%%%%%
\section{Proof of Theorem \ref{main thm-3}}%\label{sec-6}
%%%%%%%%%%%%%%%%%%%
%%%%%%%%%%%%%%%%%%%
This section Theorem \ref{main thm-3} is proved. The proof of sufficiency is divided into two lemmas. The first one is a generalization of \cite[Proposition 7.12]{BHK}. Both bounded and unbounded domains are in our considerations and our proof is also different.

\bdefe
Let  $C\geq 1$ be a constant. We say that a domain $G\subsetneq \mathbb{R}^n$ is
\begin{enumerate}  \item {\it $C$-LLC$_1$}, if for all $x\in G$ and $r>0$, then every pair of points in $B(x,r)\cap G$ can be joined by a curve in $B(x,C r)\cap G$.

  \item {\it $C$-LLC$_2$}, if for all $x\in G$ and $r>0$, then every pair of points in $G\setminus \overline{B}(x,r)$ can be joined by a curve in  $G\setminus \overline{B}(x,r/C)$.
 \item {\it $C$-LLC}, if it is both $C$-LLC$_1$ and $C$-LLC$_2$.
\end{enumerate}
\edefe

\br\label{re-1}
A $C$-LLC domain is $c$-bounded turning for $c=2C$.
\er

\begin{lem}\label{main thm-2}
Let $G\subsetneq \mathbb{R}^n$ be a $\delta$-hyperbolic, $C_0$-LLC domain. Then $G$ is $A$-uniform with $A=A(\delta, C_0, n)$.
\end{lem}
\bpf
To show that $G$ is uniform, it suffices to prove that $G$ is a diameter uniform domain (see \cite{Ma80}). For any pair of points $x, y\in G$, there exists a quasihyperbolic geodesic $\gamma=[x, y]_k$ joining $x$ and $y$. We see from Theorem \ref{main thm-1} and Remark \ref{re-1} that $G$ is $C_1$-bounded turning with $C_1=2 C_0$, and
\be\label{eq5-8}\diam(\gamma)\leq C_{po}\inf_{\alpha}\diam(\alpha)\leq 2C_0 C_{po}|x-y|,\ee
where $C_{po}$ is the constant of Theorem \ref{main thm-1}, and the infimum is taken over all curves $\alpha$ in $G$ with end points $x$ and $y$.

By \eqref{eq5-8}, we only need to check that $\gamma$ satisfies the diameter double cone condition. For all $z\in \gamma$, because the subcurves $\gamma[x, z]$ and $\gamma[z, y]$ are also quasihyperbolic geodesics, it follows from Theorem \ref{main thm-1} and Remark \ref{re-1} that
\be\label{eq5-8b}
\diam(\gamma[x, z])\leq 2C_0 C_{po}|x-z|\,\,\,\mbox{and}\,\,\,\diam(\gamma[z, y])\leq 2C_0 C_{po}|y-z|.
\ee
Denote
\be\label{eq5-8a}
\min\{\diam(\gamma[x, z]), \diam(\gamma[z, y])\}=Bd(z).
\ee
It remains to find an upper bound for the constant $B$. By \eqref{eq5-8b} and \eqref{eq5-8a}, we have
$$Bd(z)\leq 2C_0 C_{po}\min\{|x-z|, |z-y|\},$$
which implies
$$x, y\in G\setminus \overline{B}\Big(z, \frac{B}{4C_0 C_{po}}d(z)\Big).$$
 Because $G$ is $C_0$-LLC$_2$, there is another curve
 $$\beta\subset G\setminus \overline{B}\Big(z, \frac{B}{4C_0^2 C_{po}}d(z)\Big)$$
  joining $x$ and $y$. This shows that
 \be\label{eq5-8c}
 \dist(z, \beta)\geq \frac{B}{4C_0^2C_{po}}d(z).
 \ee

 By Theorem \ref{BB03-thm0.1}, $G$ satisfies the separation property with the constant $C_{sp}=C_{sp}(n, \delta)$. Thus we have
$$\beta\cap B_{\ell_G}(z, C_{sp}d(z))\neq \emptyset,$$
where $\ell_G$ denotes the length metric of $G$ associated to the $n$-Euclidean metric $|\cdot|$. Hence, it follows from \eqref{eq5-8c} that
$$\frac{B}{4C_0^2C_{po}}d(z)\leq \dist(z, \beta)\leq \dist_{\ell_G}(z, \beta)\leq C_{sp}d(z),$$
and therefore,
$$B\leq 4C_0^2C_{po}C_{sp}.$$
\epf

Next, we prove the following result by using ideas from \cite[Proposition 7.13]{BHK}. For the completeness of our proof, we give the details.

\begin{lem}\label{main lem-1}
Let $G\subsetneq \mathbb{R}^n$ be a $\delta$-hyperbolic domain. Suppose that $G$ is unbounded and there is a natural $\theta$-quasisymmetric identification
$$\varphi:\, (\partial G, |\cdot|)\to (\partial_\infty G\setminus\{\xi_0\}, d_{b, \varepsilon_0}),$$
where $\partial_\infty G$ is the Gromov boundary of $(G, k)$, $\varphi(\infty)=\xi_0$, $b\in \mathcal{B}(\xi_0)$ is a Busemann function, and $d_{b, \varepsilon_0}$ is a Hamenst\"adt metric based at $b$ with parameter $\varepsilon_0=\varepsilon_0(\delta) >0$. Then $G$ is $C$-LLC with $C=C(\theta, \delta, n)$.
\end{lem}
\bpf
Let $(G, d_\varepsilon)$ be the conformal deformation of $(G, k)$ induced by \eqref{b-0} for some fixed constant $\varepsilon=\varepsilon(\delta)>0$. It follows from Lemma \ref{main lem-2}(\ref{lem2-6}) that the identity map
$$\psi:\, (\partial_\infty G\setminus \{\xi_0\}, d_{b, \varepsilon_0})\to (\partial G_\varepsilon, d_\varepsilon)$$
is $\eta_0$-quasisymmetric with $\eta_0=\eta_0(\delta)$. Because the identity map
$$\varphi:\,(\partial G, |\cdot|)\to (\partial_\infty G\setminus\{\xi_0\}, d_{b, \varepsilon_0})$$
 is $\theta$-quasisymmetric, by \cite[Theorem 2.2]{TV}, we know that the identity map $$\psi\circ \varphi:\,(\partial G, |\cdot|)\to (\partial G_\varepsilon, d_\varepsilon)$$
  is $\theta_1$-quasisymmetric, where $\theta_1(t)=\eta_0(\theta(t))$.

First, we show that $G$ is $C$-LLC$_1$. For any point $x\in G$, $r>0$, and for all pair of points $y, z\,\in B(x, r)\cap G$, we assume that $y$ and $z$ can not be joined by a curve in $G\cap B(x, C r)$. We only need to find an upper bound for the constant $C$. Without loss of generality, we may assume that $C\geq 10$.

Let $\gamma_1$ be any curve in $G$ joining $y$ and $z$. By the assumption, we may find $y_1\in S(x, r)\cap \gamma_1$ such that $\gamma_1[y, y_1]\subset \overline{B}(x, r)$, and find a point $y_2\in S(x, \sqrt{C}r)\cap \gamma_1$ such that $\gamma_1[y, y_2]\subset \overline{B}(x, \sqrt{C}r)$. Similarly, there is a point $z_1\in S(x, r)\cap \gamma_1$ such that $\gamma_1[z, z_1]\subset \overline{B}(x, r)$ and there is $z_2\in S(x, \sqrt{C} r)\cap \gamma_1$ such that $\gamma_1[z, z_2]\subset \overline{B}(x, \sqrt{C} r)$. The connectedness of Euclidean sphere ensures that there is a curve $\alpha_1$ joining $y_1$ and $z_1$ in $S(x, r)$ and a curve $\beta_1$ joining $y_2$ and $z_2$ in $S(x, \sqrt{C}r)$ on $\partial G\cap \alpha_1$, respectively.

Because $y$ and $z$ can not be connected by any curve in $G\cap B(x, Cr)$, $\alpha_1\cap \partial G\neq \emptyset$, and $\beta_1\cap \partial G\neq \emptyset$, we may choose $u_1,v_1\,\in \alpha_1\cap \partial G$ to be the first point and the last point on $\alpha_1\cap \partial G$, when traveling from $y_1$ to $z_1$, and $u_2, v_2\,\in \beta_1\cap \partial G$ to be the first point and the last point on $\beta_1\cap \partial G$, when traveling from $y_2$ to $z_2$, respectively.
Note that we may have $u_1=v_1$ and $u_2=v_2$ (see Figure \ref{f3}).

\begin{figure}[htbp]
\begin{center}
\includegraphics{fig2.1}
\end{center}
\caption{} \label{f3}
\end{figure}

Denote
$$E_1=\alpha_1(u_1, y_1)\cup\gamma_1[y_1, y_2]\cup \beta_1(y_2, u_2),$$
and
$$F_1=\alpha_1(v_1, z_1)\cup \gamma_1[z_1, z_2]\cup \beta_1(z_2, v_2).$$
Thus $E_1\cup F_1\subset G\cap \overline{B}(x, \sqrt{C}r)$. Because $y$ and $z$ can not be connected in $G\cap B(x, Cr)$, for any curve $\sigma_1$ joining $E_1$ and $F_1$ in $G$, there is a subcurve $\sigma_{1,1}\subset \sigma_1$ between $S(x, \sqrt{C}r)$ and $S(x, Cr)$. Hence, by \cite[(2.8)]{HK}, we obtain
\be\label{eq5-3a}
\mbox{mod}_n (E_1, F_1; G)\leq \mbox{mod}_n (B(x, \sqrt{C}r), G\setminus B(x, Cr); G).
\ee

Because $u_1$, $u_2$, $v_1$, and $v_2\in \partial G$, the quasisymmetry of the map $(\partial G, |\cdot|)\to (\partial_\varepsilon G, d_\varepsilon)$ implies that

\beqq
\Delta_\varepsilon(E_1, F_1)&=&\frac{\dist_\varepsilon(E_1, F_1)}{\min\{\diam_\varepsilon(E_1),\diam_\varepsilon(F_1)\}} \\
&\leq& \max\bigg\{\frac{d_\varepsilon(u_1, v_2)}{d_\varepsilon(u_1, u_2)}, \frac{d_\varepsilon(u_1, v_2)}{d_\varepsilon(v_1, v_2)}\bigg\}\\
&\leq& \max\bigg\{\theta_1\bigg(\frac{|u_1-v_2|}{|u_1-u_2|}\bigg), \theta_1\bigg(\frac{|u_1-v_2|}{|v_1-v_2|}\bigg)\bigg\}\\
&\leq& \theta_1\Big(\frac{2\sqrt{C}}{\sqrt{C}-1}\Big)\leq \theta_1(4).
\eeqq
Then it follows from Lemma \ref{main lem-2}(\ref{lem2-4}) that
\be\label{eq5-4}
\mbox{mod}_n(E_1, F_1; G_\varepsilon)\geq \phi(\Delta_\varepsilon(E_1, F_1))\geq \phi(\theta_1(4)),
\ee
where $\phi=\phi(n, \delta)$ is the control function of Lemma \ref{main lem-2}(\ref{lem2-4}).

On the other hand, it follows from Lemma \ref{main lem-2} (\ref{lem2-3b}) that
\be\label{eq5-5}
\mbox{mod}_n(E_1, F_1; G_\varepsilon)=\mbox{mod}_n(E_1, F_1; G).
\ee
It follows from \eqref{eq5-3a} and the standard estimate of the modulus \cite[7.5]{Vai75} that
\beqq
\mbox{mod}_n(E_1, F_1; G)\leq \mbox{mod}_n(B(x, r), G\setminus B(x, \sqrt{C}r); G)\leq \omega_{n-1}(\log \sqrt{C})^{1-n},
\eeqq
where $\omega_{n-1}$ is the $(n-1)$-dimensional surface area of the unit sphere $\mathbb{S}^{n-1}$.
Therefore, this together with \eqref{eq5-4} and \eqref{eq5-5}, shows that
$$C\leq C(\omega_{n-1}, \phi, \theta_1)=C(\theta, \delta, n).$$

Next we show that $G$ is $C$-LLC$_2$. The argument is similar to the first part, but for completeness we show the details. For any point $x\in G$, $r>0$ and for all pair of points $y, z\,\in G\setminus \overline{B}(x, r)$, we assume that $y$ and $z$ can not be joined by a curve in $G\setminus \overline{B}(x, r/C)$. We only need to find an upper bound for the constant $C$. Without loss of generality, we may assume that $C\geq 10$.

Let $\gamma_2$ be any curve joining $y$ and $z$ in $G$. By the assumption we may find points $y_3\in S(x, r)\cap \gamma_2$ and $y_4\in S(x, r/\sqrt{C})\cap \gamma_2$ such that $\gamma_2[y, y_3]\subset G\setminus B(x, r)$ and $\gamma_2[y, y_4]\subset G\setminus B(x, r/\sqrt{C})$. Similarly, we may find points $z_3\in S(x, r)\cap \gamma_2$ and $z_4\in S(x, r/\sqrt{C})\cap \gamma_2$ such that $\gamma_2[z, z_3]\subset G\setminus B(x, r)$ and $\gamma_2[z, z_4]\subset G\setminus B(x, r/\sqrt{C})$.

The connectedness of Euclidean sphere implies that we can find curves $\alpha_2$ joining $y_3$ and $z_3$ in $S(x, r)$ and $\beta_2$ joining $y_4$ and $z_4$ in $S(x, r/\sqrt{C})$, respectively. Because $y$ and $z$ can not be joined within $G\setminus \overline{B}(x, r/C)$, we get $\alpha_2\cap \partial G\neq \emptyset$ and $\beta_2\cap \partial G\neq \emptyset.$ Hence, we may choose $u_3$ and $v_3$ be the first point and the last point on $\alpha_2\cap \partial G$, when traveling from $y_3$ to $z_3$, respectively. Similarly, we may choose $u_4$ and $v_4$ to be the first point and the last point on $\beta_2\cap \partial G$, when traveling from $y_4$ to $z_4$, respectively. Note that we may have $u_3=v_3$ and $u_4=v_4$. See Figure \ref{f4}.

\begin{figure}[htbp]
\begin{center}
\includegraphics{fig4.1}
\end{center}
\caption{} \label{f4}
\end{figure}

Denote
$$E_2=\alpha_2(u_3, y_3)\cup \gamma_2[y_3, y_4]\cup \beta_2(y_4, u_4)$$
and
$$F_2=\alpha_2(v_3, z_3)\cup \gamma_2[z_3, z_4]\cup \beta_2(z_4, v_4).$$
Thus $E_2\cup F_2\subset G\setminus B(x, r/\sqrt{C})$. Because $y$ and $z$ can not be connected in $G\setminus \overline{B}(x, r/C)$, thus for any curve $\sigma_2$ joining $E_2$ and $F_2$ in $G$, there is a subcurve $\sigma_{2, 2}\subset \sigma_2$ between $S(x, r/C)$ and $S(x, r/\sqrt{C})$. Hence, by \cite[(2.8)]{HK}, we obtain that
$$\mbox{mod}_n (E_2, F_2; G)\leq \mbox{mod}_n\big(B(x, r/C), G\setminus{B(x, r/\sqrt{C})}; G\big).$$
Because $u_3$, $u_4$, $v_3$, and $v_4\in \partial G$, the quasisymmetry of $(\partial G, |\cdot|)\to (\partial_\varepsilon G, d_\varepsilon)$ implies that
\beqq
\Delta_\varepsilon(E_2, F_2)&=&\frac{\dist_\varepsilon(E_2, F_2)}{\min\{\diam_\varepsilon(E_2),\diam_\varepsilon(F_2)\}} \\
&\leq& \max\Big\{\frac{d_\varepsilon(u_3, v_4)}{d_\varepsilon(u_3, u_4)}, \frac{d_\varepsilon(u_3, v_4)}{d_\varepsilon(v_3, v_4)}\Big\}\\
&\leq& \max\bigg\{\theta_1\bigg(\frac{|u_3-v_4|}{|u_3-u_4|}\bigg), \theta_1\bigg(\frac{|u_3-v_4|}{|v_3-v_4|}\bigg)\bigg\}\\
&\leq& \theta_1\Big(\frac{2\sqrt{C}}{\sqrt{C}-1}\Big)\leq \theta_1(4).
\eeqq
Then it follows from Lemma \ref{main lem-2}(\ref{lem2-4}) that
\be\label{eq5-6}
\mbox{mod}_n(E_2, F_2; G_\varepsilon)\geq \phi(\Delta_\varepsilon(E_2, F_2))\geq \phi(\theta_1(4)).
\ee

On the other hand, it follows from Lemma \ref{main lem-2}(\ref{lem2-3b}) that
\be\label{eq5-7}
\mbox{mod}_n(E_2, F_2; G_\varepsilon)=\mbox{mod}_n(E_2, F_2; G).
\ee
By the standard estimate of modulus, we obtain from \eqref{eq2-0} in Definition \ref{de-1} that
\beqq
\mbox{mod}_n(E_2, F_2; G)\leq \omega_{n-1}(\log \sqrt{C})^{1-n},
\eeqq
where $\omega_{n-1}$ is the area measure of the unit sphere $\mathbb{S}^{n-1}$.
This together with \eqref{eq5-6} and \eqref{eq5-7} implies that
$$C\leq C(\omega_{n-1}, \phi, \theta)=C(\theta, \delta, n).$$
Hence the proof of Lemma \ref{main lem-1} is complete.
\epf

 Now we are ready to prove Theorem \ref{main thm-3}. The following result is also needed.

\begin{lem}$($\cite[Lemma 3.14]{BHK}$)$\label{BHK lem 3.14}
Let $G\subsetneq \mathbb{R}^n$ be an $A$-uniform domain, let $x,$ $y$ and $z\in G$ satisfing $|x-y|\geq 2|x-z|$, and let $\gamma$ and $\alpha$ be the quasihyperbolic geodesics joining $x$ to $y$ and $z$, respectively. Let $w\in \gamma$ be such that  $\ell(\gamma[x, w])=|x-z|$. Then
$$\dist_k(y, \alpha)-C_A\leq k(y, w)\leq \dist_k(y, \alpha)+C_A,$$
where the constant $C_A$ depends on $A$.
\end{lem}
\br\label{re-2}
In \cite[Lemma 3.14]{BHK}, $G$ is assumed to be bounded. However, we find from a carefully checking its proof that this requirement is not necessary in our case.
\er

\subsection*{Proof of Theorem \ref{main thm-3}}
The sufficiency follows from Lemma \ref{main thm-2} and Lemma \ref{main lem-1}. It suffices to show the necessity part. Assume that $G$ is an unbounded $A$-uniform domain. Then \cite[Theorem 3.6]{BHK} shows that $G$ is $\delta$-hyperbolic with $\delta=\delta(A)$. Moreover, it follows from \cite[Proposition 3.12]{BHK} that there is a natural identification $\varphi:\,\partial G\to \partial_\infty G\setminus\{\xi_0\}$ such that $\varphi(\infty)=\xi_0$ for some $\xi_0\in \partial_\infty G$. So we only need to prove that $\varphi$ is $\eta$-quasisymmetric, where $\partial_\infty G\setminus \{\xi_0\}$ is equipped with a Hamenst\"adt metric $d_{b, \varepsilon_0}$ for a Busemann function $b=b_{o,\xi_0}\in \mathcal{B}(\xi_0)$, $o\in G$, and $\varepsilon_0=\varepsilon_0(\delta)$.

Fix three distinct points $x, y, z\in \partial G$. Denote $|x-y|=t|x-z|$, we need to find a homeomorphism $\eta:\,[0, \infty)\to [0,\infty)$ such that
\be\label{eq7-0}
d_{b, \varepsilon_0}(\varphi(x), \varphi(y))\leq \eta(t)d_{b,\varepsilon_0}(\varphi(x), \varphi(z)).
\ee
Choose sequences $\{x_n\}\to x$, $\{y_n\}\to y$, and $\{z_n\}\to z$ in the Euclidean metric $|\cdot|$ as $n\to \infty$, respectively. Also choose a sequence $\{u_n\}\subset G$ with $|u_n-o|\to \infty$ as $n\to \infty$. We see from the arguments in \cite[Proposition 3.12]{BHK} that $\{x_n\}\in \varphi(x)$, $\{y_n\}\in \varphi(y)$, $\{z_n\}\in \varphi(z)$ and $\{u_n\}\in \varphi(\infty)=\xi_0$ are Gromov sequences.
By \eqref{zz0.1} and Lemma \ref{z0}, we have
\beqq
(x_n|z_n)_b-(x_n|y_n)_b&\leq& (x_n|z_n)_o-(\xi_0|z_n)_o-(x_n|y_n)_o+(y_n|\xi_0)_o+20\delta\\
&\leq& (x_n|z_n)_o-(u_n|z_n)_o-(x_n|y_n)_o+(y_n|u_n)_o+22\delta\\
&=&(x_n|z_n)_{u_n}-(x_n|y_n)_{u_n}+22\delta.
\eeqq
This together with Lemma \ref{BuSc-Lem3.2.4} shows that
\beqq
(\varphi(x)|\varphi(z))_b-(\varphi(x)|\varphi(y))_b&\leq& \liminf\limits_{n\to \infty}\big((x_n|z_n)_b-(x_n|y_n)_b\big)+44\delta\\
&\leq& \liminf\limits_{n\to \infty} \big((x_n|z_n)_{u_n}-(x_n|y_n)_{u_n}\big)+66\delta.
\eeqq

Now by \eqref{eq-b}, we obtain
\beq\label{eq7-3}
T:= \frac{d_{b, \varepsilon_0}(\varphi(x), \varphi(y))}{d_{b, \varepsilon_0}(\varphi(x), \varphi(z))}\leq 2e^{66\delta \varepsilon_0}\liminf\limits_{n\to \infty}e^{\varepsilon_0(x_n|z_n)_{u_n}-\varepsilon_0(x_n|y_n)_{u_n}}.
\eeq
Fix a sufficiently large positive integer $n$. Set
$$s_n=(x_n|z_n)_{u_n}-(x_n|y_n)_{u_n}\,\,\,\mbox{and}\,\,\,t_n=
\frac{|x_n-y_n|}{|x_n-z_n|}.$$
Because $|u_n-o|\to \infty$ and $|x_n-x|\to 0$ as $n\to \infty$, respectively, we may assume without loss of generality that for all $n$,
\be\label{eq7-4}
\max\{|x_n-y_n|, |x_n-z_n|\}\leq \frac{1}{2}|x_n-u_n|.
\ee
Choose quasihyperbolic geodesics $\alpha_n$, $\beta_n$, and $\gamma_n$ connecting $x_n$ to the points $u_n$, $y_n$, and $z_n$, respectively. By \eqref{eq7-4}, there are two points $v_n$, $w_n\in \alpha_n$ such that
\be\label{eq7-5}
\ell(\alpha_n[x_n, v_n])=|x_n-y_n|\,\,\,\mbox{and}\,\,\,\ell(\alpha_n[x_n, w_n])=|x_n-z_n|,
\ee
respectively. We see from \eqref{eq7-4}, \eqref{eq7-5}, and Lemma \ref{BHK lem 3.14} that
\be\label{eq7-5a}
\dist_k(u_n, \beta_n)-C_A\leq k(u_n, v_n)\leq \dist_k(u_n, \beta_n)+C_A
\ee
\be\label{eq7-5b}
\dist_k(u_n, \gamma_n)-C_A\leq k(u_n, w_n)\leq \dist_k(u_n, \gamma_n)+C_A,
\ee
where $C_A$ is the constant of Lemma \ref{BHK lem 3.14} with $C_A=C(A)$.

Because it follows from \cite[Theorem 2.10]{BHK} that $\alpha_n$ is $C_{A}'$-uniform with $C_{A}'=C_{A}'(A)$, we obtain from \eqref{eq7-5} that
\be\label{eq7-6}
d(v_n)\geq \frac{1}{C_{A}'}\ell(\alpha_n[x_n, v_n])= \frac{1}{C_{A}'}|x_n-y_n|
\ee
and
\be\label{eq7-7}
d(w_n)\geq \frac{1}{C_{A}'}\ell(\alpha_n[x_n, w_n])= \frac{1}{C_{A}'}|x_n-z_n|.
\ee
In the following, we divide the arguments into two cases: $t\geq 1$ and $0<t<1$.

Suppose first that $t\geq 1$. Without loss of generality, we may assume that $t_n\geq 1/2$ for all $n$, because $t_n\to t$ as $n\to \infty$. Thus we have
\be\label{eq7-8}
|x_n-y_n|\geq \frac{1}{2}|x_n-z_n|.
\ee
From \eqref{eq7-8}, \eqref{eq7-6}, and \eqref{eq7-7}, it follows that
\be\label{eq7-9}
\min\{d(v_n), d(w_n)\}\geq \frac{1}{2C_{A}'}|x_n-z_n|.
\ee
We see from \eqref{eq7-5} and \eqref{eq7-8} that
\be\label{eq7-10}
|v_n-w_n|\leq |v_n-x_n|+|x_n-w_n|\leq |x_n-y_n|+|x_n-z_n|\leq 3|x_n-y_n|.
\ee
Thus it follows that
\beqq
s_n&=&(x_n|z_n)_{u_n}-(x_n|y_n)_{u_n}\\
&\leq& \dist_k(u_n, \gamma_n)-\dist_k(u_n, \beta_n)+2\delta \,\,\,\,\quad\quad\quad\quad \quad(\mbox{By\,\,Lemma}\,\,\,\ref{2.33})\\
&\leq& k(u_n, w_n)-k(u_n, v_n)+2C_{A}+2\delta\,\,\,\,\quad\quad\quad\quad\quad (\mbox{By}\,\,\, \eqref{eq7-5a}\,\,\,\mbox{and}\,\,\,\eqref{eq7-5b})\\
&\leq& k(v_n, w_n)+2C_A+2\delta\\
&\leq&4A^2\log\bigg(1+\frac{|v_n-w_n|}{\min\{d(v_n), d(w_n)\}}\bigg)+2C_A+2\delta\,\,\,(\mbox{By\,\,Theorem \,\,\ref{BHK 2.4}})\\
&\leq& 4A^2\log(1+6C_{A}'t_n)+2C_A+2\delta.\,\,\,\,\,\,\quad\quad\quad\quad\quad\quad
(\mbox{By}\,\,\,\eqref{eq7-9})
\eeqq
This together with \eqref{eq7-3}, shows that
$$T\leq \liminf_{n\to \infty} 2e^{68\delta\varepsilon_0+2C_A\varepsilon_0}e^{4 A^2 \varepsilon_0\log(1+16C_{A}'t_n)}=Ct^{4 A^2\varepsilon_0}.$$

Next, we consider the case that $0<t<1$. We may assume that $t_n\leq 1$ for large enough $n$. This guarantees that $w_n\in \alpha_n[v_n, u_n]$ by the choices of $u_n$ and $w_n$ in \eqref{eq7-6} and \eqref{eq7-7}, respectively.  We have
\beq\label{eq7-11}
s_n&=&(x_n|z_n)_{u_n}-(x_n|y_n)_{u_n}\\
&\leq& \dist_k(u_n, \gamma_n)-\dist_k(u_n, \beta_n)+2\delta \,\,\,\,\quad \quad(\mbox{By\,\,Lemma}\,\,\,\ref{2.33})\nonumber\\
&\leq& k(u_n, w_n)-k(u_n, v_n)+2C_A+2\delta\,\,\,\,\quad\quad (\mbox{By}\,\,\, \eqref{eq7-5a}\,\,\,\mbox{and}\,\,\,\eqref{eq7-5b})\nonumber\\
&=& -k(v_n, w_n)+2C_A+2\delta, \nonumber
\eeq
where the last equality follows from the fact that $w_n$ lies on the quasihyperbolic geodesic $\alpha_n[v_n, u_n]$.

On the other hand, by \eqref{eq7-5}, \eqref{eq7-6}, and \eqref{eq7-7}, we get
\be\label{eq7-12}
d(v_n)\leq |v_n-x|\leq |v_n-x_n|+|x_n-x|\leq |x_n-y_n|+|x_n-x|.
\ee
Because $|x_n-x|\to 0$ as $n\to \infty$, we may assume without loss of generality that  $|x_n-x|\leq |x_n-y_n|$. It follows from \eqref{eq7-7}, \eqref{eq7-12}, and Theorem \ref{BHK 2.4} that
$$-k(v_n, w_n)\leq \log\frac{d(v_n)}{d(w_n)}\leq \log\frac{2C_{A}'|x_n-y_n|}{|x_n-z_n|}=\log(2C_{A}' t_n).$$
Hence, by \eqref{eq7-3} and \eqref{eq7-11}, we obtain
$$T\leq \liminf_{n\to \infty} 2e^{68\delta\varepsilon_0+2C_A\varepsilon_0 }e^{\varepsilon_0 \log (2C_{A}'t_n)}\leq Ct^{\varepsilon_0}.$$
\qed

%%%%%%%%%%%%%%%%%%%
%%%%%%%%%%%%%%%%%%%
\bigskip
%
%{\bf Acknowledgement.} The authors would like to thank Professors M. Vuorinen and Y. Li for many valuable comments and suggestions. The authors are indebted to the referee for the valuable suggestions.

%%%%%%%%%%%%%%%%%%%

\end{document}